# CENTRAL LIMIT THEOREMS FOR ITERATED RANDOM LIPSCHITZ MAPPINGS

By Hubert Hennion and Loïc Hervé

*Université de Rennes 1 and Institut National des Sciences Appliquées de Rennes*

Let $M$ be a noncompact metric space in which every closed ball is compact, and let $G$ be a semigroup of Lipschitz mappings of $M$. Denote by $(Y_n)_{n \geq 1}$ a sequence of independent $G$-valued, identically distributed random variables (r.v.'s), and by $Z$ an $M$-valued r.v. which is independent of the r.v. $Y_n$, $n \geq 1$. We consider the Markov chain $(Z_n)_{n \geq 0}$ with state space $M$ which is defined recursively by $Z_0 = Z$ and $Z_{n+1} = Y_{n+1} Z_n$ for $n \geq 0$. Let $\xi$ be a real-valued function on $G \times M$. The aim of this paper is to prove central limit theorems for the sequence of r.v.'s $(\xi(Y_n, Z_{n-1}))_{n \geq 1}$. The main hypothesis is a condition of contraction in the mean for the action on $M$ of the mappings $Y_n$; we use a spectral method based on a quasi-compactness property of the transition probability of the chain mentioned above, and on a special perturbation theorem.

**1. Introduction.** Let $M$ be a noncompact metric space in which every closed ball is compact, endowed with its Borel $\sigma$-field $\mathcal{M}$. We denote by $G$ a semigroup of Lipschitz mappings of $M$ and by $\mathcal{G}$ a $\sigma$-field on $G$. We assume that the action of $G$ on $M$ is measurable; that is, the map $j$ defined by $j(g, y) = gy$ is measurable from $(G \times M, \mathcal{G} \otimes \mathcal{M})$ to $(M, \mathcal{M})$.

Let $\pi$ be a probability distribution on $G$, and let $(Y_n)_{n \geq 1}$ be a sequence of independent $G$-valued random variables (r.v.'s) identically distributed according to $\pi$, defined on a probability space $(\Omega, \mathcal{F}, \mathbb{P})$. The iterated random mappings $R_n$, $n \geq 0$, are defined by

$$R_0 = Id_M, \qquad R_n = Y_n \cdots Y_1, \ n \geq 1.$$

Let $Z$ be an $M$-valued r.v. which is independent of the r.v.'s $Y_n$, $n \geq 1$. The sequence $(Z_n)_{n \geq 0}$ defined by

$$Z_n = R_n Z, \qquad n \geq 0,$$

---

is a Markov chain on $M$ which is defined recursively by

$$Z_0 = Z, \qquad Z_{n+1} = j(Y_{n+1}, Z_n) = Y_{n+1}Z_n, \ n \geq 0.$$

Observe that we get here the general Lipschitz iterative model on $M$, which has been considered by many authors; see Duflo (1997) and Diaconis and Freedman (1999) to get an overview of the subject. Consider particularly the case where $M$ is the linear space $\mathbb{R}^q$. The generalized linear autoregressive model is obtained when $G$ is the semigroup of affine mappings of $M$. Replace in the preceding the linear part of the action by that of a fixed Lipschitz mapping $f$ of $M$. An element $g$ of the semigroup $G$ is now defined by a vector $b_g \in M$, and it acts on $M$ according to the formula $gx = f(x) + b_g$. In this context the probability distribution $\pi$ on $G$ is simply defined by a distribution on $M$; thus we get the Lipschitz functional autoregressive model.

Now let $\xi$ be a real-valued measurable function on $G \times M$. The aim of this paper is to establish a central limit theorem with a rate of convergence and a local central limit theorem for the sequence of r.v.'s

$$(\xi(Y_n, Z_{n-1}))_{n \geq 1}.$$

The interest of considering a function $\xi$ of the couple $(g, x) \in G \times M$ rather than a function only depending on $x$ appears, for example, in the study of random matrices products.

From the stochastic viewpoint, the context may be described as the study of the sequence of r.v.'s obtained by composing the function $\xi$ and the Markov chain $(X_n)_{n \geq 0}$ with state space $G \times M$ defined by

$$X_0 = (Id_M, Z), \qquad X_n = (Y_n, Z_{n-1}), \qquad n \geq 1.$$

The main hypothesis will be a condition of contraction in the mean of the action on $M$ of the elements of $G$ under the probability distribution $\pi$. This property enables us to make use of a refinement of the spectral method. Recall that the spectral method was initiated by Nagaev (1957), and then used and improved by many authors. It is fully described in Hennion and Hervé (2001), where references are given. The spectral method is based on a quasi-compactness property of the transition probability $Q$ of the chain $(X_n)_{n \geq 0}$, and on a perturbation theorem ensuring that, for small $|t|$, the Fourier kernels $Q(t)$ associated with $Q$ and $\xi$ have spectral properties similar to those of $Q$. In the present setting, the use of the standard perturbation theory for operators leads to assume moments of exponential type (cf. Milhaud and Raugi (1989) and Hennion and Hervé (2001), Chapter X, Section 3). The main feature of this paper is the use of a perturbation theorem of Keller and Liverani (1999) which is adapted to operators verifying a Doeblin–Fortet inequality. By means of this theorem, we get the desired limit theorems under moments of polynomial types.

Notice that there are several methods to cope with central limit theorems for a function of a Markov chain; most known are regeneration and splitting, use of central limit theorems for martingale increments and Lindeberg techniques. As will be discussed later, when applied to the present context, some of these methods can give a central limit theorem under hypotheses which are weaker than ours; however, it seems that these methods have not yet been developed so far as to get the central limit theorem with a rate of convergence and the local central limit theorem of this paper. See Section 3 for more details.

**2. Statements of results.** For $g \in G$, we set

$$c(g) = \sup \left\{ \frac{d(gx, gy)}{d(x, y)} : x, y \in M, \ x \neq y \right\},$$

by assumption $c(g) < +\infty$.

For $n \in \mathbb{N}^*$, we denote by $\pi^{*n}$ the distribution of $R_n$. We choose a fixed point $x_0$ in $M$. For $\eta \geq 1$ and $n \in \mathbb{N}^*$, we define the integrals:

$$\mathcal{M}_\eta = \int_G (1 + c(g) + d(gx_0, x_0))^\eta \, d\pi(g),$$

$$\mathcal{M}'_\eta = \int_G c(g) \, (1 + c(g) + d(gx_0, x_0))^{\eta-1} \, d\pi(g),$$

$$\mathcal{C}^{(n)}_\eta = \int_G c(g) \max\{c(g), 1\}^{\eta-1} \, d\pi^{*n}(g).$$

Notice that, since $c(\cdot)$ is submultiplicative, $\mathcal{M}'_\eta < +\infty$ implies $\mathcal{C}^{(n)}_\eta < +\infty$.

The statements below will appeal, on the one hand to the moment conditions $\mathcal{M}_\eta < +\infty$ and $\mathcal{M}'_{\eta'} < +\infty$, on the other hand, to the average contractivity condition $\mathcal{C}^{(n)}_{\eta'} < 1$, for a suitable choice of $\eta, \eta' \geq 1$.

We consider a real-valued measurable function $\xi$ on $G \times M$ satisfying:

CONDITION RS. *There exist $r, s \in \mathbb{R}_+$ and measurable nonnegative functions $R, S$ on $G$ such that, for all $x, y \in M$ and $g \in G$,*

$$|\xi(g, x)| \leq R(g)(1 + d(x, x_0))^r,$$

$$|\xi(g, x) - \xi(g, y)| \leq S(g) \, d(x, y)(1 + d(x, x_0) + d(y, x_0))^s.$$

Observe that, if the second condition in Condition RS holds, then the first one is also valid with $r = s + 1$ and $R(g) = |\xi(g, x_0)| + S(g)$. However, it is worth noticing that this condition may be verified for a smaller exponent $r$; this is the case, for example, when $\xi$ is bounded. This remark also shows that, without a significant loss of generality, we could add to Condition RS the inequality $r \leq s + 1$; yet, we notice that, when $r$ increases, $R(g)$ decreases.

The case $s = 0$ and $r = 1$ corresponds to functions $\xi$ such that $\xi(g, \cdot)$ is Lipschitz for all $g \in G$. At last, notice that, if $\alpha \in ]0, 1]$, then $d(\cdot, \cdot)^\alpha$ is a distance on $M$; consequently, Condition RS involves the case of functions $\xi$ such that $\xi(g, \cdot)$ is locally $\alpha$-Hölder for all $g \in G$.

As in the Introduction, we denote by $Z$ a r.v. in $M$ defined on $(\Omega, \mathcal{F}, P)$, and independent of the r.v. $Y_n$, $n \geq 1$. We set

$$S_n^Z = \sum_{k=1}^n \xi(Y_k, Z_{k-1}), \qquad n \geq 1.$$

We now state central limit theorems for the sequence $(S_n^Z)_n$; more precise results concerning the behavior of the sequence $(R_n Z, S_n^Z)_n$ are given in Section 9.

A preliminary to all these statements is the existence of a probability distribution on $M$ which is preserved by the action of $\pi$. More precisely, the action on $M$ of the sequence of random mappings $(R_n)_{n \geq 0}$ defines a Markov chain: for $y_0 \in M$, the sequence $(R_n y_0)_{n \geq 0}$ is Markov with state space $M$, initial distribution $\delta_{y_0}$, and transition probability $P$ defined by

$$y \in M, \ B \in \mathcal{M}, \qquad P(y, B) = \int_G \mathbb{1}_B(gy) \, d\pi(g).$$

THEOREM I (Invariant probability measure).   *Assume that there exist $\gamma \geq 0$ and an integer $n_0 \geq 1$ such that $\mathcal{M}_{\gamma+1} < +\infty$ and $\mathcal{C}_{\gamma+1}^{(n_0)} < 1$.*

*Then there exists on $(M, \mathcal{M})$ a unique $P$-invariant probability distribution $\nu$. Moreover, we have*

$$\int_M d(x, x_0)^{\gamma+1} \, d\nu(x) < +\infty,$$

*and the geometric ergodicity holds in the Prohorov distance $d_{\mathrm{P}}$. Namely, there exist positive real numbers $C$ and $\kappa_0 < 1$, such that, for any probability distribution $\mu$ on $M$ satisfying $\mu(d(\cdot, x_0)) < +\infty$, and all $n \geq 1$,*

$$d_P(\mu P^n, \nu) \leq C \kappa_0^{n/2}.$$

It must be noted that such an ergodicity result holds under much weaker hypotheses; see the survey of Diaconis and Freedman (1999) and a recent result in Bhattacharya and Majumdar (2004). In fact, the above statement is just the one which fits the general hypotheses of the paper.

In the sequel our hypotheses will involve a parameter $\gamma_0 > 0$ and:

CONDITIONS $\mathcal{H}(\gamma_0)$.   For

$$\mathcal{M}_{\gamma_0 + 1} < +\infty, \qquad \mathcal{M}'_{2\gamma_0 + 1} < +\infty,$$

there exists $n_0 \in \mathbb{N}^*$ such that $\mathcal{C}_{2\gamma_0 + 1}^{(n_0)} < 1$.

Since $\mathcal{C}_{\gamma_0+1}^{(n_0)} \leq \mathcal{C}_{2\gamma_0+1}^{(n_0)}$, if the above conditions hold, then the $P$-invariant distribution $\nu$, whose existence is ensured by Theorem I, is such that $\nu(d(\cdot, x_0)^{\gamma_0+1}) < +\infty$; consequently, if the number $r$ and the function $R$ in Condition RS verify $r \leq \frac{\gamma_0}{2} + \frac{1}{2}$ and $\int_G R(g)^2 \, d\pi(g) < +\infty$, we have

$$\int_M \int_G \xi(g, x)^2 \, d\pi(g) \, d\nu(x) < +\infty.$$

From now on we shall assume that

$$m = \int_M \int_G \xi(g, x) \, d\pi(g) \, d\nu(x) = 0.$$

This causes no loss of generality since it is always possible to replace $\xi$ by $\xi - m$. Otherwise, we shall keep in mind that, if $Z$ has the $P$-invariant distribution $\nu$, then we have $\mathbb{E}[d(Z, x_0)^{\gamma_0+1}] < +\infty$. However, unless otherwise stated, in the sequel $Z$ is not supposed to be $\nu$-distributed.

At last, we define for $g \in G$,

$$\tilde{\delta}(g) = 1 + c(g) + d(gx_0, x_0),$$

and for $\tau > 0$ and positive real valued measurable functions $U$, $V$ on $G$, we set

$$\mathcal{J}^\tau(U, V) = \int_G U(g)c(g)\tilde{\delta}(g)^{2\tau} \, d\pi(g) + \int_G V(g)\tilde{\delta}(g)^{\tau+1} \, d\pi(g),$$

or more shortly $\mathcal{J}^\tau(U, V) = \pi(Uc\,\tilde{\delta}^{2\tau}) + \pi(V\,\tilde{\delta}^{\tau+1})$.

**THEOREM A** (Central limit). *Assume $\mathcal{H}(\gamma_0)$ with $\gamma_0 > r + \max\{r, s+1\}$ and that*

$$\int_G R^2 \, d\pi < +\infty, \qquad \mathcal{J}^{\gamma_0-r}(R, R+S) < +\infty.$$

*Then there exists a real number $\sigma^2 \geq 0$ such that, under the condition $\mathbb{E}[d(Z, x_0)^{\gamma_0+1}] < +\infty$, the sequence $(\frac{S_n^Z}{\sqrt{n}})_{n \geq 1}$ converges in distribution to a $\mathcal{N}(0, \sigma^2)$-distributed r.v.*

As already mentioned, this statement is not the best known one; using our spectral method, it is a stage to the two following results.

**THEOREM B** (Central limit with a rate of convergence). *Assume $\mathcal{H}(\gamma_0)$ with $\gamma_0 > 3r + \max\{r, s+1\}$ and that*

$$\int_G R^3 \, d\pi < +\infty, \qquad \mathcal{J}^{\gamma_0-r}(R, R+S) + \mathcal{J}^{\gamma_0-2r}(R^2, (R+S)R) < +\infty.$$

*Then, if $\sigma^2 > 0$, there exists a constant $C$ such that, when $Z$ verifies $\mathbb{E}[d(Z, x_0)^{\gamma_0+1}] < +\infty$, we have, for all $n \geq 1$,*

$$\sup_{u \in \mathbb{R}} |\mathbb{P}[S_n^Z \leq u\sigma\sqrt{n}] - \mathcal{N}(0, 1)(]-\infty, u])| \leq C \, \frac{1 + \mathbb{E}[d(Z, x_0)^{\gamma_0+1}]}{\sqrt{n}}.$$

We denote by $\mathcal{L}$ the Lebesgue measure on $\mathbb{R}$. Furthermore, a complex-valued function on $M$ is said to be locally Lipschitz if it is Lipschitz on every compact subset of $M$.

THEOREM C (Local central limit).    *Assume that the conditions of Theorem* A *hold, and that $\xi$ verifies the nonarithmeticity condition: there is no $t \in \mathbb{R}$, $t \neq 0$, no $\lambda \in \mathbb{C}$, $|\lambda| = 1$, no bounded locally Lipschitz function $w$ on $M$ with nonzero constant modulus on the support $\Sigma_\nu$ of $\nu$, such that we have, for all $x \in \Sigma_\nu$ and all $n \geq 1$,*

$$e^{itS_n^x} w(R_n x) = \lambda^n w(x), \qquad \mathbb{P}\text{-}a.s.$$

*Then, if $\sigma^2 > 0$, and if $Z$ is such that $\mathbb{E}[d(Z, x_0)^{\gamma_0 + 1}] < +\infty$, we have, for every continuous function $h$ on $\mathbb{R}$ such that $\lim_{|u| \to +\infty} u^2 h(u) = 0$,*

$$\lim_n \sigma \sqrt{2\pi n} \mathbb{E}[h(S_n^Z)] = \mathcal{L}(h).$$

We end with a result which gives a criterion for $\sigma^2 > 0$ and defines $\sigma^2$ asymptotically.

THEOREM S.    *Assume $\mathcal{H}(\gamma_0)$ with $\gamma_0 > 2r + s + 1$ and that*

$$\mathcal{J}^{\gamma_0 - r}(R, R + S) + \mathcal{J}^{\gamma_0 - 2r}(R^2, (R + S)R) < +\infty.$$

(i) *If $\sigma^2 = 0$, then there exists a real-valued locally Lipschitz function $\tilde{\xi}_1$ on $M$ satisfying $\nu(\tilde{\xi}_1^2) < +\infty$, and such that we have, with $Z$ distributed according to $\nu$,*

$$\xi(Y_1, Z) = \tilde{\xi}_1(Z) - \tilde{\xi}_1(Y_1 Z), \qquad \mathbb{P}\text{-}a.s.$$

(ii) *If the distribution of $Z$ verifies $\mathbb{E}[d(Z, x_0)^{\gamma_0 + 1}] < +\infty$, then*

$$\sigma^2 = \lim_n \frac{1}{n} \mathbb{E}[(S_n^Z)^2].$$

It will be seen later on (Theorems C′ and S′) that the functions $w$ and $\tilde{\xi}_1$ in the two last statements must not be merely locally Lipschitz; they must belong to certain spaces to be defined in the sequel.

In the following section, we show how these theorems apply to some cases of interest. This being done, the rest of the paper is devoted to the proofs; the reader will find in Section 4.2 a brief outline of the subsequent work.

## 3. Applications.

3.1. *Sequences of type* $(\chi(Z_n))_n$. Let $\chi$ be a real-valued locally Lipschitz function on $M$, and suppose that there exist $C, s \in \mathbb{R}_+$ such that, for all $x, y \in M$,

$$|\chi(x) - \chi(y)| \leq C d(x,y)(1 + d(x,x_0) + d(y,x_0))^s.$$

Using martingale methods, it is proved that the central limit theorem for $(\chi(Z_n))_n$ holds for any initial distribution under the moment condition $\int_G d(gx_0, x_0)^{4(s+1)} \, d\pi(g) < +\infty$ and the contraction property $\int_G c(g)^{4(s+1)} \, d\pi(g) < 1$; see Duflo (1997). By means of similar techniques, it is established in Benda (1998) that, when $s = 0$, the same result is valid under the weaker hypotheses $\int_G c(g)^2 \, d\pi(g) < 1$ and $\int_G d(gx_0, x_0)^2 \, d\pi(g) < +\infty$. Considering the stationary chain with initial probability $\nu$, Wu and Woodroofe (2000) have established a central limit theorem for functions $\chi$ which are not Lipschitz and not even continuous.

Let us now apply the results of the preceding section: we set $\xi(g, x) = \chi(x)$. The moment hypotheses of Theorem A are the same as those of Theorem C, so that it can be seen from Theorem 3.1 that they are stronger than the ones previously stated. However, to our knowledge, Theorems B and C are new. They can be stated as follows. Recall that $\tilde{\delta}(g) = 1 + c(g) + d(gx_0, x_0)$.

THEOREM 3.1. *Suppose that there exist* $\varepsilon > 0$ *and integers* $n_0 \geq 1$, $k \geq 0$ *such that*

$$\pi(\tilde{\delta}^{k(s+1)+1+\varepsilon/2} + c\tilde{\delta}^{2k(s+1)+\varepsilon}) < +\infty,$$

*and*

$$\pi^{*n_0}(c \max\{1, c\}^{2k(s+1)+\varepsilon}) < 1,$$

*and assume that* $\nu(\chi) = 0$, *where* $\nu$ *is the* $P$-*invariant probability measure.*

*If* $k$ *takes the values* 4 *and* 3, *respectively, then the assertions of Theorems* B *and* S, *respectively, apply to*

$$S_n^Z = \sum_{k=1}^{n} \chi(Z_{k-1}).$$

*Moreover, if* $\chi$ *is nonarithmetic and if the above integral conditions are satisfied for* $k = 2$, *then the assertion of Theorem* C *holds.*

PROOF. The function $\xi$ on $G \times M$ defined by $\xi(g, y) = \chi(y)$ verifies Condition RS with the exponents $r = s + 1$ and $s$ associated with constant functions $R$ and $S$. These have moments of all orders. Consequently the moment conditions of Theorems B, S, and C reduce to $\mathcal{H}(\gamma_0)$ with $\gamma_0 = k(s+1) + \frac{\varepsilon}{2}$; this gives the desired results. □

Let us point out that Pollicott (2001) has stated a central limit theorem with a rate of convergence and a large deviations theorem in the case where the support of the probability measure $\pi$ is finite. However, this study is based on the assertion without proof that, on a suitable space of Lipschitz functions, the Fourier kernels $P(t)$ (see Section 4) are analytic perturbed operators of $P$. Also notice that, if it is proved that the stationary chain with initial probability distribution $\nu$ is strongly mixing and Harris recurrent, then we can apply Bolthausen (1982) to obtain a central limit theorem with an $n^{-1/2}$ rate of convergence. However, on one hand, this requires some additional hypotheses on $\pi$ [see Meyn and Tweedie (1993), page 140, for a sufficient condition in the context of the following section]; on the other hand, this only covers the stationary case.

3.2. *Generalized autoregressive processes.* Denote by $G$ the semigroup of all affine mappings of $M = \mathbb{R}^q$, $q \geq 1$. An element $g \in G$ is identified with a couple $(a(g), b(g))$, where $a(g)$ is an endomorphism of $\mathbb{R}^q$ and $b(g)$ is a vector in $\mathbb{R}^q$. For $y \in M$, we set $gy = a(g)y + b(g)$. The associated generalized autoregressive process $(Z_n)_{n \geq 0}$ is then defined by

$$Z_0 = Z, \qquad Z_{n+1} = a(Y_{n+1})Z_n + b(Y_{n+1}), \qquad n \geq 0.$$

Let $\xi$ be a function from $G \times \mathbb{R}^q$ to $\mathbb{R}$, and suppose that there exist a norm $\|\cdot\|$ on $\mathbb{R}^q$, $\alpha \in ]0,1]$, $r, s \in \mathbb{R}_+$ and nonnegative measurable functions $R$ and $S$ on $G$ such that, for all $g \in G$ and $x, y \in \mathbb{R}^q$, we have

$$|\xi(g, x)| \leq R(g)(1 + \|x\|)^{\alpha r},$$

$$|\xi(g, x) - \xi(g, y)| \leq S(g)\|x - y\|^{\alpha}(1 + \|x\| + \|y\|)^{\alpha s}.$$

For instance, these properties hold with $\alpha = 1$ when $\xi$ is a polynomial function of the entries of the matrix representing $a(g)$ and of the coordinates of the vectors $b(g)$ and $x$.

Let us consider the distance $d$ defined on $\mathbb{R}^q$ by $d(x, y) = \|x - y\|^{\alpha}$, and choose $x_0 = 0 \in \mathbb{R}^q$. We have $c(g) = \|a(g)\|^{\alpha}$ and $d(gx_0, x_0) = \|b(g)\|^{\alpha}$. Then the statements B, C, S apply straightforwardly. To compare with former results, let us rewrite Theorem B. Let $\tilde{\delta}(g) = (1 + \|a(g)\| + \|b(g)\|)^{\alpha}$, then

THEOREM 3.2. *The hypotheses in the central limit theorem with a rate of convergence (Theorem* B*) are satisfied if there exist* $\gamma_0 > 3r + \max\{r, s+1\}$ *and* $n_0 \in \mathbb{N}^*$ *such that*

$$\pi(\tilde{\delta}^{\gamma_0 + 1} + \|a\|^{\alpha}\tilde{\delta}^{2\gamma_0}) < +\infty \quad \text{and} \quad \pi^{*n_0}(\|a\|^{\alpha}\max\{1, \|a\|\}^{2\gamma_0\alpha}) < 1,$$

*and when the functions* $R(\cdot)$ *and* $S(\cdot)$ *satisfy the moment conditions*

$$\int_G R^3 \, d\pi < +\infty, \qquad \mathcal{J}^{\gamma_0 - r}(R, R + S) + \mathcal{J}^{\gamma_0 - 2r}(R^2, (R + S)R) < +\infty.$$

In this context, convergence rates in the central limit theorem have already been established by Milhaud and Raugi (1989) and by Cuny (2004). The spectral method used in Milhaud and Raugi (1989) is, in substance, similar to the one developed here, but it appeals to the standard perturbation theorem. For this reason (see Section 6.1), the following conditions on $a(\cdot)$ and $b(\cdot)$ are required: $\|a(\cdot)\| < 1$ $\pi$-p.s, and there exist real numbers $\rho > 0$ and $\beta \in \,]0, 1]$ such that we have the exponential moment condition

$$\int_G e^{\rho \|b(g)\|^\beta} (R(g) + S(g))^5 (1 - \|a(g)\|^\beta)^{-5\alpha/\beta(1+\max\{r,s+1\})} \, d\pi(g) < +\infty.$$

The hypotheses on both $a(\cdot)$ and $b(\cdot)$ are significantly less restrictive in Theorem 3.2. The study in Cuny (2004) is based on martingale methods. The contraction condition is the same as in Milhaud and Raugi (1989), and it is supposed that, for all $\ell \in \mathbb{N}$, $\int \|b(g)\|^\ell \, d\pi(g) < +\infty$. Under these conditions, for functions $\xi$ which are not necessarily Hölder of the variable $x$, it is proved that the rate of convergence in the central limit theorem is $n^{-p}$ for every $p < \frac{1}{2}$.

3.3. *Products of positive random matrices.* Let $G$ be the semigroup of $q \times q$ matrices with nonnegative entries which are allowable, namely, every row and every column contains a strictly positive element, and denote by $G^\circ$ the ideal of $G$ composed of matrices with strictly positive entries.

For $g \in G$ and $w \in \mathbb{R}^q$, we denote by $g(w)$ the image of $w$ under $g$; the cone

$$C = \{w : w = (w_1, \ldots, w_q) \in \mathbb{R}^q, w_k > 0, \, k = 1, \ldots, q\}$$

is invariant under all $g \in G$. Define $M$ to be the intersection of the hyperplane $\{w : w \in \mathbb{R}^q, \sum_{k=1}^q w_k = 1\}$ of $\mathbb{R}^q$ with $C$.

The linear space $\mathbb{R}^q$ is endowed with the norm $\|\cdot\|$ defined by

$$w = (w_1, \ldots, w_q) \in \mathbb{R}^q, \qquad \|w\| = \sum_{k=1}^q |w_k|,$$

and for each $g \in G$, we set

$$\|g\| = \sup\{\|g(y)\| : y \in M\}, \qquad v(g) = \inf\{\|g(y)\| : y \in M\}.$$

The semigroup $G$ being equipped with its Borel $\sigma$-field $\mathcal{G}$, we consider a probability distribution $\pi$ on $G$ for which there exists an integer $n_0$ such that the support of the r.v. $R_{n_0}$ contains a matrix of $G^\circ$. Denote by $g^*$ the adjoint of $g$. It is shown in Hennion (1997) that, if

$$\int_G (|\ln \|g^*\|\,| + |\ln v(g^*)|)^2 \, d\pi(g) < +\infty,$$

then there exists $\gamma_1 \in \mathbb{R}$ such that, for $y \in M$, the sequence $(\frac{1}{\sqrt{n}}(\ln \|R_n(y)\| - n\gamma_1))_{n \geq 1}$ converges to the $\mathcal{N}(0, \sigma^2)$ distribution; moreover, the case $\sigma^2 = 0$ is investigated. Using Theorems B and C, it is possible to state a central limit theorem with a rate of convergence and a local central limit theorem. Notice that similar theorems have already been given in Hennion and Hervé [(2001), Section X.5], but under more restrictive moment hypotheses.

To see how this case enters the present frame, we set, for $g \in G$ and $y \in M$,

$$gy = \frac{g(y)}{\|g(y)\|}, \qquad a(g, y) = \ln \|g(y)\|.$$

It is easy to check that the first formula defines an action of $G$ on $M$, while the function defined by the second one verifies the property of additive cocycle associated with this action:

$$a(gg', y) = a(g, g'y) + a(g', y), \qquad g, g' \in G, \ y \in M.$$

Consequently, setting $\xi(g, y) = a(g, y) - \gamma_1$, for $(g, y) \in G \times M$, we can write, for $y \in M$,

$$\ln \|R_n(y)\| - n\gamma_1 = \sum_{k=1}^{n} \xi(Y_k, R_{k-1}y).$$

Furthermore, when $M$ is endowed with a suitable metric $d_H$ called the Hilbert metric [see Bapat and Raghavan (1997)], every $g \in G$ is Lipschitz with constant $c(g) \leq 1$, and we have $c(g) < 1$ if and only if all entries of $g$ are strictly positive. Therefore, if the support of $R_{n_0}$ contains such a matrix $g$, we have $C_\eta^{(n_0)} = \int_G c(g) \, d\pi^{*n_0}(g) < 1$ for all $\eta \geq 1$.

For $\eta \geq 0$, set

$$\mathcal{L}^\eta = \int_G (|\ln \|g\|| + |\ln v(g)| + |\ln v(g^*)|)^\eta \, d\pi(g).$$

THEOREM 3.3. *Suppose that there exists an integer $n_0$ such that the support of the r.v. $R_{n_0}$ contains a matrix of $G^\circ$, and let $\varepsilon > 0$.*

(i) *Assume $\mathcal{L}^{4+\varepsilon} < +\infty$; then, if $\sigma^2 > 0$, there exists a nonnegative constant $C$ such that, in case the r.v. $Z$ of $M$ verifies $\mathbb{E}[d_H(Z, x_0)^{2+\varepsilon/2}] < +\infty$, we have, for all $n \geq 1$,*

$$\sup_{u \in \mathbb{R}} |\mathbb{P}[\ln \|R_n(Z)\| - n\gamma_1 \leq u\sigma\sqrt{n}] - \mathcal{N}(0,1)(]-\infty, u])|$$

$$\leq C \frac{1 + \mathbb{E}[d_H(Z, x_0)^{2+\varepsilon/2}]}{\sqrt{n}}.$$

(ii) *Assume $\mathcal{L}^{3+\varepsilon} < +\infty$, $\sigma^2 > 0$, and that the support of $R_{n_0}$ contains two matrices $g_1, g_2 \in G^\circ$ whose spectral radii $\rho_1$, $\rho_2$ verify $\ln \frac{\rho_2}{\rho_1} \notin \mathbb{Q}$. Then, if $\mathbb{E}[d_H(Z, x_0)^{2+\varepsilon/2}] < +\infty$, we have, for any real valued continuous function $h$ on $\mathbb{R}$ such that $\lim_{|u| \to +\infty} u^2 h(u) = 0$,*

$$\lim_n \sigma \sqrt{2\pi n}\, \mathbb{E}[h(\ln \|R_n(Z)\| - n\gamma_1)] = \mathcal{L}(h).$$

PROOF. First notice that the number $\|g^*\|$ associated to any endomorphism $g$ of $\mathbb{R}^q$ defines a new norm which is equivalent to the one already considered. Consequently, there exists a constant $C$ such that, for every $g \in G$, we have $|\ln \|g^*\|\,| \le C + |\ln \|g\|\,|$.

We denote by $(e_k)_{k=1}^q$ and by $\langle \cdot, \cdot \rangle$ the canonical basis and scalar product on $\mathbb{R}^q$. If $y, y' \in M$, we set

$$m_H(y, y') = \min\left\{ \frac{\langle y, e_k \rangle}{\langle y', e_k \rangle} : k = 1, \ldots, q \right\},$$

$$d_H(y, y') = -\ln(m_H(y, y') m_H(y', y)),$$

$d_H$ is the Hilbert distance on $M$ [see Bapat and Raghavan (1997)]. The space $(M, d_H)$ is not compact, but each closed ball in it is compact. Set $x_0 = (1/q, \ldots, 1/q) \in M$; we have

$$d(gx_0, x_0) = \ln \frac{\max_i \langle gx_0, e_i \rangle}{\min_j \langle gx_0, e_j \rangle} = \ln \frac{\max_i \|g^* e_i\|}{\min_j \|g^* e_j\|}$$

$$= \ln \frac{\|g^*\|}{v(g^*)} \le C + |\ln \|g\|\,| + |\ln v(g^*)\,|.$$

The function $\xi(g, \cdot)$ is bounded by $|\gamma_1| + |\ln \|g\|| + |\ln v(g)|$. From $\|gy\| \ge m_H(y, y') \|gy'\|$, we deduce that $S(g) = 1$ and $s = 0$ [see Hennion (1997), Lemma 5.3]. Therefore, Condition RS is verified with $R(g) = 2(|\ln \|g\|| + |\ln v(g)|)$, $r = 0$, and $S(g) = 1$, $s = 0$.

The above estimations prove that the required moment conditions of Theorems B and C hold if we have, respectively, $\mathcal{L}^{4+\varepsilon} < +\infty$ and $\mathcal{L}^{3+\varepsilon} < +\infty$.

It remains to prove that the additional hypothesis in (ii) implies the nonarithmeticity of $\xi$. Let $k = 1, 2$. It follows from the Perron–Frobenius theorem that $\rho_k > 0$ and that, for all $\ell \ge 1$, we have $g_k^\ell = \rho_k^\ell (p_k + h_k^\ell)$, where $p_k \in G^\circ$ and the endomorphism $h_k$ of $\mathbb{R}^q$ has a spectral radius $< 1$. Consequently, for any $x \in M$, we have $\ln \|g_k^\ell x\| = \ell \ln \rho_k + r_{k,\ell}(x)$, with $\lim_\ell r_{k,\ell}(x) = \ln \|p_k(x)\|$. Suppose that there exist $t \in \mathbb{R}$, $t \ne 0$, $\lambda \in \mathbb{C}$, $|\lambda| = 1$, and a bounded locally Lipschitz function $w$ on $M$ which has a nonzero constant modulus on the support $\Sigma_\nu$ of $\nu$, and such that we have, for all $x \in \Sigma_\nu$ and all $n \ge 1$,

$$\lambda^n w(x) = e^{itS_n^x} w(R_n x) = e^{it(\ln \|R_n(x)\| - n\gamma_1)} w(R_n x), \qquad \mathbb{P}\text{-a.s.}$$

From the continuity of the functions used in the two members, we deduce that, for any $\ell \geq 1$ and $x \in \Sigma_\mu$, we have

$$e^{it(\ln \|g_k^\ell(x)\| - n_0 \ell \gamma_1)} w(g_k^\ell x) = \lambda^{n_0 \ell} w(x).$$

It follows that

$$e^{it\ell \ln(\rho_2/\rho_1)} = \frac{w(g_1^\ell x)}{w(g_2^\ell x)} e^{it(r_{1,n_0\ell}(x) - r_{2,n_0\ell}(x))}.$$

The second member converges when $\ell \to +\infty$, while the countable set of complex numbers defined by the first one is dense in $\{z : z \in \mathbb{C}, |z| = 1\}$. This contradiction completes the proof.  $\square$

## 4. Preliminaries.

4.1. *$P$-invariant probability measure ( proof of Theorem* I). Here the hypotheses are those of Theorem I: $\mathcal{M}_{\gamma+1} < +\infty$, $\mathcal{C}_{\gamma+1}^{(n_0)} < 1 (\gamma \geq 0, n_0 \in \mathbb{N}^*)$. For $\lambda \in ]0, 1]$, $x \in M$, and $g \in G$, we set

$p_\lambda(x) = 1 + \lambda d(x, x_0),$

$\delta_\lambda(g) = \max\{c(g), 1\} + \lambda d(gx_0, x_0)$  and  $\tilde{\delta}(g) = 1 + c(g) + d(gx_0, x_0).$

LEMMA 4.1.  *We have for all $g \in G$ and $0 \leq \lambda \leq 1$*

$$\sup_{x \in M} \frac{p_\lambda(gx)}{p_\lambda(x)} \leq \delta_\lambda(g) \leq \tilde{\delta}(g).$$

*The functions $c(\cdot)$ and $\tilde{\delta}(\cdot)$ are submultiplicative.*

PROOF.  Let $x \in M$ and $g \in G$; then

$$\frac{p_\lambda(gx)}{p_\lambda(x)} = \frac{1 + \lambda d(gx, x_0) - \lambda d(gx_0, x_0)}{1 + \lambda d(x, x_0)} + \frac{\lambda d(gx_0, x_0)}{1 + \lambda d(x, x_0)}$$

$$\leq \frac{1 + \lambda d(gx, gx_0)}{1 + \lambda d(x, x_0)} + \lambda d(gx_0, x_0) \leq \frac{1 + \lambda c(g) d(x, x_0)}{1 + \lambda d(x, x_0)} + \lambda d(gx_0, x_0)$$

$$\leq \max\{1, c(g)\} + \lambda d(gx_0, x_0).$$

The fact that $c(\cdot)$ is submultiplicative is obvious. Finaly for $h, g \in G$, we get

$$\tilde{\delta}(hg) \leq 1 + c(h)c(g) + [d(hgx_0, x_0) - d(hx_0, x_0)] + d(hx_0, x_0)$$

$$\leq 1 + c(h)c(g) + c(h)d(gx_0, x_0) + d(hx_0, x_0) \leq \tilde{\delta}(h)\tilde{\delta}(g). \qquad \square$$

Recall that, for $n \in \mathbb{N}^*$, $\pi^{*n}$ denotes the law of $R_n$.

LEMMA 4.2. *Let* $\phi_\lambda(x) = d(x, x_0)p_\lambda(x)^\gamma$ *for* $\lambda \in [0, 1]$*. Then:*

(a) *For all* $n \geq 1$ *and* $x \in M$*, we have* $P^n\phi_\lambda(x) < +\infty$*.*
(b) *For* $\lambda_0 \in\ ]0, 1]$ *small enough, we have* $\int_G c(g)\delta_{\lambda_0}(g)^\gamma\, d\pi^{*n_0}(g) < 1$*.*
(c) *There exist constants* $\varepsilon \in\ ]0, 1[$*,* $C \in \mathbb{R}_+$ *such that*

$$P^{n_0}\phi_{\lambda_0} \leq C + \varepsilon\phi_{\lambda_0}.$$

PROOF. (a) For $n \geq 1$ and $x \in M$, we have

$$P^n\phi_\lambda(x) = \int_G d(gx, x_0)p_\lambda(gx)^\gamma\, d\pi^{*n}(g)$$

$$\leq \int_G d(gx_0, x_0)p_\lambda(gx)^\gamma\, d\pi^{*n}(g)$$

$$\qquad + \int_G [d(gx, x_0) - d(gx_0, x_0)]p_\lambda(gx)^\gamma\, d\pi^{*n}(g)$$

$$\leq p_\lambda(x)^\gamma \int_G d(gx_0, x_0)\delta_\lambda(g)^\gamma\, d\pi^{*n}(g)$$

$$\qquad + d(x, x_0)p_\lambda(x)^\gamma \int_G c(g)\delta_\lambda(g)^\gamma\, d\pi^{*n}(g).$$

The functions in the two integrals above are dominated by $\tilde{\delta}(\cdot)^{\gamma+1}$. Since this function is submultiplicative and $\pi$-integrable, Fubini's theorem ensures that these integrals are finite. Thus $P^n\phi_\lambda(x) < +\infty$.

Since $c\delta_{\lambda_0}^\gamma \leq \tilde{\delta}^{\gamma+1}$ and $\tilde{\delta}$ is submultiplicative, assertion (b) is a direct consequence of hypotheses and Lebesgue's theorem.

Set $\varepsilon' = \int_G c(g)\delta_{\lambda_0}(g)^\gamma\, d\pi^{*n_0}(g)$. From the above inequality applied with $\lambda = \lambda_0$ and $n = n_0$, there exists a constant $D_0$ such that

$$P^{n_0}\phi_{\lambda_0} \leq D_0 p_{\lambda_0}^\gamma + \varepsilon'\phi_{\lambda_0}.$$

Using continuity and $\lim_{d(x, x_0) \to +\infty} \frac{p_{\lambda_0}(x)^\gamma}{\phi_{\lambda_0}(x)} = 0$, we see that there exists a constant $C$ such that $D_0 p_{\lambda_0}^\gamma \leq C + \frac{1-\varepsilon'}{2}\phi_{\lambda_0}$. Hence $P^{n_0}\phi_{\lambda_0} \leq C + \varepsilon\phi_{\lambda_0}$ with $\varepsilon = \frac{1+\varepsilon'}{2}$. $\square$

Now let us prove Theorem I. For convenience we set $\phi = \phi_{\lambda_0}$. By induction and Lemma 4.2, we obtain, for every $q \geq 1$, $P^{qn_0}\phi \leq \varepsilon^q\phi + C(1 + \varepsilon + \cdots + \varepsilon^{q-1}) \leq \phi + \frac{C}{1-\varepsilon}$. Let $n \in \mathbb{N}^*$. Writing $n = qn_0 + r$ with $r \in \{0, \ldots, n_0 - 1\}$ and setting $E = \max\{P^k\phi(x_0),\ k = 0, \ldots, n_0 - 1\}$, we get $P^n\phi(x_0) \leq E + \frac{C}{1-\varepsilon}$. Therefore, the sequence $(P^n\phi(x_0))_n$ is bounded by a constant, say $K$. For $n \geq 1$, let $\nu_n$ be the probability measure on $(M, \mathcal{M})$ defined by

$$B \in \mathcal{M}, \qquad \nu_n(B) = \frac{1}{n}\sum_{k=0}^{n-1}(P^k\mathbb{1}_B)(x_0).$$

Observe that, for each $n \geq 1$, we have $\nu_n(\phi) \leq K$. Since $\lim_{d(x,x_0) \to +\infty} \phi(x) = \infty$, the subset $[\phi \leq \alpha]$ is compact for each $\alpha > 0$. The Markov inequality implies that, for all $n \geq 1$, we have $\nu_n([\phi > \alpha]) \leq \frac{\nu_n(\phi)}{\alpha} \leq \frac{K}{\alpha}$, so that the sequence $(\nu_n)_n$ is tight. Therefore, we can select a subsequence $(\nu_{n_k})_k$ converging to a probability measure $\nu$. It is clear that $\nu$ is $P$-invariant.

For $p \in \mathbb{N}^*$, set $\phi_p(\cdot) = \min(\phi(\cdot), p)$. For $k \geq 0$ and $p \geq 0$, we have $\nu_{n_k}(\phi_p) \leq \nu_{n_k}(\phi) \leq K$; consequently, for all $p \geq 0$, $\lim_k \nu_{n_k}(\phi_p) \leq K$. The monotone convergence theorem gives $\nu(\phi) < +\infty$, that is, $\nu(d(\cdot, x_0)^{\gamma+1}) < +\infty$.

Now let us prove that $\nu$ is the unique $P$-invariant probability distribution. First observe that, since $\mathbb{E}[\ln c(R_{n_0})] \leq \ln \mathbb{E}[c(R_{n_0})] = \ln \mathcal{C}_1^{(n_0)} < 0$, the law of large numbers asserts that $\limsup_q c(R_{qn_0})^{1/q} \leq \lim_q (\prod_{\ell=1}^{q} c(Y_{\ell n_0} \cdots Y_{(\ell-1)n_0+1}))^{1/q} < 1$ on a set $\Omega_1$ such that $\mathbb{P}(\Omega_1) = 1$. For $x, y \in M$ and $q \geq 1$, we can write $d(R_{qn_0}x, R_{qn_0}y) \leq c(R_{qn_0})\, d(x,y)$, so that $\lim_q d(R_{qn_0}x, R_{qn_0}y) = 0$ on $\Omega_1$. Let $\nu'$ be a $P$-invariant probability distribution on $M$. For each bounded continuous function $f$ on $M$, we have

$$\nu'(f) - \nu(f) = \int_M E[f(R_{qn_0}x) - f(R_{qn_0}y)]\, d\nu'(x)\, d\nu(y);$$

passing to the limit, we get $\nu'(f) - \nu(f) = 0$. We conclude that $\nu' = \nu$.

It remains to establish the geometric ergodicity in the Prohorov distance $d_{\mathrm{P}}$. Let $f$ be a bounded uniformly Lipschitz function on $M$. Then, for all $x, y \in M$ and $n \geq 1$, we have

$$|P^n f(x) - P^n f(y)| \leq \int_G |f(gx) - f(gy)|\, d\pi^{*n}(g)$$

$$\leq m_0(f) d(x,y) \int_G \frac{d(gx, gy)}{d(x,y)}\, d\pi^{*n}(g)$$

$$\leq m_0(f)\, d(x,y) \mathcal{C}_1^{(n)},$$

where $m_0(f) = \sup\{\frac{|f(x)-f(y)|}{d(x,y)}, x, y \in M, \ x \neq y\}$. Let $\mu$ be the law of $Z$, and assume that $\phi_0(\cdot) = d(\cdot, x_0)$ is $\mu$-integrable. By integrating the previous inequality with respect to both $d\nu(x)$ and $d\mu(y)$, it follows that $|\nu(f) - \mu P^n(f)| \leq \mathcal{C}_1^{(n)} m_0(f)(\nu(\phi_0) + \mu(\phi_0))$. This bound proves that $\nu - \mu P^n$ is a continuous linear functional on the space of all bounded uniformly Lipschitz functions on $M$ endowed with its canonical norm. Moreover, we have $\|\nu - \mu P^n\| \leq C' \mathcal{C}_1^{(n)}$ with $C' = \nu(\phi_0) + \mu(\phi_0)$. Writing $n = qn_0 + r$ with $r \in \{0, \ldots, n_0-1\}$ and using the fact that $c(\cdot)$ is submultiplicative, we easily see that $\mathcal{C}_1^{(n)} \leq C''(\mathcal{C}_1^{(n_0)})^{n/n_0}$. Since $d_p(\nu, \mu P^n) \leq 2\|\nu - \mu P^n\|^{1/2}$ [see Dudley (1989)], the last assertion of Theorem I follows with $\kappa_0 = (\mathcal{C}_1^{(n_0)})^{1/n_0}$.

4.2. *Outlines of the method.* As mentioned in the Introduction, the main idea of this work consists in applying the method described in Hennion and Hervé (2001) to the function $\xi$ and to the Markov chain $(X_n)_{n \geq 0}$ with the state space $G \times M$ and the transition probability $Q$ defined by

$$(g, y) \in G \times M, \ B \in \mathcal{G} \otimes \mathcal{M}, \qquad Q((g, y), B) = \int_G \mathbb{1}_B(h, gy) \, d\pi(h).$$

However, we observed in Chapter X of Hennion and Hervé (2001), devoted to Lipschitz kernels, that, because of the special form of $Q$, the essential part of the study can be performed with the help of the transition probability $P$ and of the Fourier kernels $P(t)$, $t \in \mathbb{R}$, associated to $P$ and $\xi$, which are defined, for any bounded measurable function $f$ on $M$, by

$$y \in M, \qquad P(t)f(y) = \int_G e^{it\xi(g,y)} f(gy) \, d\pi(g).$$

This is due to the fact that, for all functions $f$ as above, we have $Q(f \circ j) = (Pf) \circ j$, where $j$ is the action of $G$ on $M$. Then, in the sequel, we shall only use the kernels $P(t)$; the next statement indicates that these kernels are sufficient for our purpose.

BASIC LEMMA. *Let $f$ be a bounded measurable function on $M$, and denote by $\mu$ the distribution of $Z$. Then we have, for $n \geq 1$, $t \in \mathbb{R}$,*

$$\mathbb{E}[f(R_n Z)e^{itS_n^Z}] = \mu(P(t)^n f).$$

PROOF. Set $S_0^Z = 0$. For $n \geq 1$, we have

$$\mathbb{E}[f(R_n Z)e^{itS_n^Z}] = \mathbb{E}[f(Y_n R_{n-1} Z)e^{it(S_{n-1}^Z + \xi(Y_n, R_{n-1} Z))}].$$

Since $(Z, Y_1, \ldots, Y_{n-1})$ and $Y_n$ are independent r.v.'s, Fubini's theorem gives

$$\mathbb{E}[f(R_n Z)e^{itS_n^Z}] = \mathbb{E}[e^{itS_{n-1}^Z} \int_G f(g R_{n-1} Z)e^{it\xi(g, R_{n-1} Z)} \, d\pi(g)]$$

$$= \mathbb{E}[e^{itS_{n-1}^Z}(P(t)f)(R_{n-1} Z)].$$

The desired formula for $n = 1$ holds because the second member equals $\mathbb{E}[P(t)f(Z)] = \mu(P(t)f)$. Suppose now that the stated formula is valid at rank $n - 1$, $n \geq 2$. Then, from the previous relation and the fact that $P(t)f$ is a bounded measurable function on $M$, we conclude that $\mathbb{E}[f(R_n Z)e^{itS_n^Z}] = \mu(P(t)^{n-1}(P(t)f))$. This completes the proof. $\square$

Theorems A, B, C, S will be direct consequences of the extensions $A', B', C', S'$ stated in Section 9. The outline of the argumentation is the following. In Section 5, we shall introduce spaces $\mathcal{B}_\gamma$, which depend on a real parameter $\gamma > 0$ and are composed of locally Lipschitz functions on $M$. Three norms, denoted by $N_{\infty,\gamma}$, $N_\gamma$ and $N_{1,\gamma}$, will be defined on $\mathcal{B}_\gamma$. It will be proved that they are equivalent, but each of them will be suited to a part of the proof. In this way, in Section 5.3, we shall see that the use of $N_\gamma$ is convenient to establish that, for suitable $\gamma$, $P$ is quasi-compact on $\mathcal{B}_\gamma$, and furthermore that the number 1 is the unique peripheral eigenvalue of $P$. In Section 6, the norms $N_{\infty,\gamma}$ will be helpful for the study of the behavior of the function $P(t)$ near $t = 0$. For this purpose, it will be worth noticing that, for $\gamma' < \gamma$, $P(t)$ may be viewed as a bounded linear map from $\mathcal{B}_{\gamma'}$ to $\mathcal{B}_\gamma$; indeed, the derivative kernels of $P(t)$, which in general do not define bounded endomorphisms of $(\mathcal{B}_\gamma, N_{\infty,\gamma})$, can be considered on the other hand as bounded linear maps from $\mathcal{B}_{\gamma'}$ to $\mathcal{B}_\gamma$ for suitable $\gamma' < \gamma$. Of course, this will be a less restrictive property because the space $\mathcal{B}_\gamma$ strictly contains $\mathcal{B}_{\gamma'}$ and is endowed with a weaker norm. In Section 7, the norm $N_{1,\gamma}$ will be an essential tool to apply a perturbation theorem due to Keller and Liverani (1999), from which it will follow that $P(t)$ are perturbed operators of $P$ for small $|t|$. The interest of this perturbation theorem is that it only requires $P(\cdot)$ to be continuous as a map taking values in the space of bounded linear map from $(\mathcal{B}_\gamma, N_{1,\gamma})$ to $(\mathcal{B}_\gamma, \nu(|\cdot|))$; this is the key point of this study (see Section 6.1). In particular, this theorem ensures that, for small $|t|$, $P(t)$ has only one dominating simple eigenvalue, $\lambda(t)$, on $\mathcal{B}_\gamma$, and we shall establish in Section 8 that the Taylor expansions for $P(t)$ at $t = 0$ obtained in Section 6 lead to expansions of the eigenelements belonging to $\lambda(t)$. Then, in Section 9, by using the previous preparation and by applying the method described in Hennion and Hervé (2001), we shall be in a position to prove limit theorems. Notice that renewal and large deviations theorems for the sequence $(S_n^Z)_{n \geq 1}$ might be derived from similar techniques.

## 5. The space $\mathcal{B}_\gamma$ and quasi-compactness of $P$.

5.1. *Conventions and notation.* From now on, we fix $\gamma_0 > 0$ and $n_0 \in \mathbb{N}^*$ such that Condition $\mathcal{H}(\gamma_0)$ holds, that is:

$$\mathcal{M}_{\gamma_0+1} = \pi(\tilde{\delta}^{\gamma_0+1}) < +\infty,$$

$$\mathcal{M}'_{2\gamma_0+1} = \pi(c\,\tilde{\delta}^{2\gamma_0}) < +\infty,$$

$$\mathcal{C}^{(n_0)}_{2\gamma_0+1} = \pi^{*n_0}(c\,\max\{c,1\}^{2\gamma_0}) < 1.$$

According to the subsequent statements, some additional conditions will be imposed on $\gamma_0$.

LEMMA 5.1. *There exists a real number $\lambda_0 \in\, ]0,1]$ such that*

$$\vartheta_0 = \int_G c(g)\,(\max\{c(g),1\} + \lambda_0\,d(gx_0,x_0))^{2\gamma_0}d\pi^{*n_0}(g) \;<\; 1.$$

PROOF. Since $c(g)(\max\{c(g), 1\} + \lambda_0 d(gx_0, x_0))^{2\gamma_0} \le c(g)\tilde{\delta}(g)^{2\gamma_0}$, and the functions $c, \tilde{\delta}$ are submultiplicative, the lemma follows from the two last conditions of $\mathcal{H}(\gamma_0)$ and Lebesgue's theorem. $\square$

Now we fix a $\lambda_0 \in \,]0, 1]$ satisfying the previous inequality.

For $x, y \in M$, $g \in G$, we set

$$p(x) = 1 + \lambda_0 d(x, x_0),$$

$$\delta(g) = \max\{c(g), 1\} + \lambda_0 d(gx_0, x_0).$$

Notice that $p \le 1 + d(\cdot, x_0) \le \frac{1}{\lambda_0} p$, and $\delta \le \tilde{\delta} \le \frac{2}{\lambda_0} \delta$. Besides, for $\gamma > 0$, let us write

$$\Delta_\gamma(x, y) = d(x, y) p(x)^\gamma p(y)^\gamma.$$

With the help of these elements, we now define the space $\mathcal{B}_\gamma$ composed of locally Lipschitz functions on $M$, and we define four equivalent norms on this space. Such spaces, introduced in Le Page (1983), have already been used by several authors in order to prove the quasi-compactness of probability kernels having a contracting property; see Milhaud and Raugi (1989) and Peigné (1993). A similar statement will be established in Section 5.3.

5.2. *Definitions of $\mathcal{B}_\gamma$ and of the norms $N_{\infty,\gamma}, N_{\infty,\gamma,\tilde{\gamma}}, N_\gamma$ and $N_{1,\gamma}$.* For $\gamma > 0$, we denote by $\mathcal{B}_\gamma$ the space of all complex-valued locally Lipschitz functions on $M$ such that

$$m_\gamma(f) = \sup\left\{\frac{|f(x) - f(y)|}{\Delta_\gamma(x, y)}, x, y \in M, x \ne y\right\} < +\infty.$$

The inequality $\Delta_\gamma(x, x_0) = d(x, x_0) p(x)^\gamma \le (1/\lambda_0) p(x)^{\gamma+1}$ ensures that, for all $f \in \mathcal{B}_\gamma$, we have $|f(x)| \le |f(x_0)| + (1/\lambda_0) m_\gamma(f) p(x)^{\gamma+1}$; thus $\sup_{x \in M} \frac{|f(x)|}{p(x)^{\gamma+1}} < +\infty$. Consequently $\mathcal{B}_\gamma$ can be equipped with the norm

$$N_{\infty,\gamma}(f) = m_\gamma(f) + |f|_\gamma,$$

where $|f|_\gamma = \sup\{\frac{|f(x)|}{p(x)^{\gamma+1}}, \ x \in M\}$.

Let $\tilde{\gamma} > \gamma$. As $p^{\tilde{\gamma}+1} \ge p^{\gamma+1}$, we have, for $f \in \mathcal{B}_\gamma$, $|f|_{\tilde{\gamma}} = \sup_{x \in M} \frac{|f(x)|}{p(x)^{\tilde{\gamma}+1}} < +\infty$; we set

$$N_{\infty,\gamma,\tilde{\gamma}}(f) = m_\gamma(f) + |f|_{\tilde{\gamma}}.$$

Since $\mathcal{M}_{\gamma_0+1} < +\infty$ and $\mathcal{C}_{\gamma_0+1}^{(n_0)} \le \mathcal{C}_{2\gamma_0+1}^{(n_0)} < 1$, the $P$-invariant probability measure, $\nu$, whose existence is ascertained by Theorem I, is such that

$$\int_M d(x, x_0)^{\gamma_0+1} \, d\nu(x) < +\infty.$$

Therefore, for every $\gamma \in\, ]0, \gamma_0]$, $\nu$ integrates $p^{\gamma+1}$, and thus integrates all the functions of $\mathcal{B}_\gamma$, so that we can define on $\mathcal{B}_\gamma$ the following norms:

$$N_\gamma(f) = m_\gamma(f) + |\nu(f)|,$$

$$N_{1,\gamma}(f) = m_\gamma(f) + \nu(|f|).$$

PROPOSITION 5.2.  *Let* $\gamma$, $0 < \gamma \leq \gamma_0$. *The four norms* $N_{\infty,\gamma}$, $N_{\infty,\gamma,\tilde{\gamma}}$, $N_\gamma$ *and* $N_{1,\gamma}$ *are equivalent on* $\mathcal{B}_\gamma$. *When equipped with one of these norms,* $\mathcal{B}_\gamma$ *is a Banach space.*

PROOF.  The fact that $(\mathcal{B}_\gamma, N_{\infty,\gamma})$ is a Banach space is well known.

(i) $N_{\infty,\gamma}$ *and* $N_{\infty,\gamma,\tilde{\gamma}}$ *are equivalent.*   Since $|f|_{\tilde{\gamma}} \leq |f|_\gamma$, we have $N_{\infty,\gamma,\tilde{\gamma}}(f) \leq N_{\infty,\gamma}(f)$. Conversely, for $x \in M$,

$$\frac{|f(x)|}{p(x)^{\gamma+1}} \leq \frac{|f(x_0)| + (1/\lambda_0)m_\gamma(f)p(x)^{\gamma+1}}{p(x)^{\gamma+1}} \leq |f(x_0)| + \frac{1}{\lambda_0}m_\gamma(f).$$

The bounds $|f(x_0)| \leq |f|_{\tilde{\gamma}}$ and $1 \leq \lambda_0^{-1}$ prove that $|f|_\gamma \leq \lambda_0^{-1} N_{\infty,\gamma,\tilde{\gamma}}(f)$; consequently, $N_{\infty,\gamma}(f) \leq (1 + \lambda_0^{-1})N_{\infty,\gamma,\tilde{\gamma}}(f)$.

To establish that $N_\gamma(\cdot)$ and $N_{\infty,\gamma}(\cdot)$ are equivalent, we proceed as in Hennion and Hervé [[2001](#), Chapter X].

LEMMA 5.3.  $(\mathcal{B}_\gamma, N_\gamma)$ *is a Banach space.*

PROOF.  Let $(f_n)_n$ be a Cauchy sequence in $(\mathcal{B}_\gamma, N_\gamma)$. Set $g_n = f_n - f_n(y_0)$, where $y_0$ is any point of $M$. We have $|g_q(x) - g_p(x)| \leq m_\gamma \times (g_q - g_p)\Delta_\gamma(x, y_0) = m_\gamma(f_q - f_p)\Delta_\gamma(x, y_0)$ because $g_n(y_0) = 0$. Hence $\nu(|g_q - g_p|) \leq \nu(\Delta_\gamma(\cdot, y_0))m_\gamma(f_q - f_p)$. Recall that $\nu(p^{\gamma+1}) < +\infty$, so that $\nu(\Delta_\gamma(\cdot, y_0)) < +\infty$.

Consequently, $(g_n)_n$ is a Cauchy sequence in the Lebesgue space $\mathbb{L}^1(\nu)$; therefore it converges in this space, and $(\nu(g_n))_n$ converges in $\mathbb{C}$. Moreover, $(\nu(f_n))_n$ converges in $\mathbb{C}$ because, by assumption, it is a Cauchy sequence. It follows that $(f_n(y_0))_n$ converges to a complex number, say $f(y_0)$, and then $(f_n)_n$ converges in $\mathbb{L}^1(\nu)$. Because $y_0$ is arbitrary, $(f_n)_n$ converges pointwise to $f$. We have $\lim_{n \to +\infty} \nu(f - f_n) = 0$. The properties $f \in \mathcal{B}_\gamma$ and $\lim_{n \to +\infty} m_\gamma(f - f_n) = 0$ are obtained by standard arguments.  □

(ii) $N_\gamma$ *and* $N_{\infty,\gamma}$ *are equivalent.*   For $f \in \mathcal{B}_\gamma$, we have $|\nu(f)| \leq \nu(|f|) \leq |f|_\gamma \nu(p^{\gamma+1})$. Thus $N_\gamma(f) \leq (1 + \nu(p^{\gamma+1}))N_{\infty,\gamma}(f)$. Since $(\mathcal{B}_\gamma, N_\gamma)$ and $(\mathcal{B}_\gamma, N_{\infty,\gamma})$ are Banach spaces, the open mapping theorem yields the claimed equivalence [see Dunford and Schwartz ([1958](#))].

(iii)  $N_{1,\gamma}$ *and* $N_{\infty,\gamma}$ *are equivalent.*    We have $|f(y)| \leq |f(x)| + m_\gamma(f)\, d(x,y) p(x)^\gamma p(y)^\gamma$ for all $x, y \in M$. By integrating this inequality with respect to the measure $\nu$, we obtain

$$|f(y)| \leq \nu(|f|) + m_\gamma(f) p(y)^\gamma \int_M d(x, x_0) p(x)^\gamma \, d\nu(x)$$

$$+ m_\gamma(f)\, d(y, x_0) p(y)^\gamma \int_M p(x)^\gamma \, d\nu(x)$$

$$\leq \nu(|f|) + 2\lambda_0^{-1} m_\gamma(f) p(y)^{\gamma+1} \nu(p^{\gamma+1}),$$

hence $N_{\infty,\gamma}(f) \leq (1 + 2\lambda_0^{-1}\nu(p^{\gamma+1})) N_{1,\gamma}(f)$.

Finally, we have $N_{1,\gamma}(f) = m_\gamma(f) + \nu(|f|) \leq m_\gamma(f) + |f|_\gamma \nu(p^{\gamma+1}) \leq (1 + \nu(p^{\gamma+1})) N_{\infty,\gamma}(f)$.  $\square$

We conclude this section by giving a statement that will be useful for the spectral study of $P(t)$.

LEMMA 5.4.    (i) *For* $0 < \gamma < \tilde{\gamma}$, *the canonical embedding from* $(\mathcal{B}_\gamma, N_{\infty,\gamma,\tilde{\gamma}})$ *into* $(\mathcal{B}_\gamma, |\cdot|_{\tilde{\gamma}})$ *is compact.*

(ii)  *For* $\gamma \in \,]0, \gamma_0]$, *the canonical embedding from* $(\mathcal{B}_\gamma, N_{1,\gamma})$ *into* $(\mathcal{B}_\gamma, \nu(|\cdot|))$ *is compact.*

PROOF.    (i) Let $(f_n)_n$ be a sequence of functions in $\mathcal{B}_\gamma$ such that $N_{\infty,\gamma,\tilde{\gamma}}(f_n) \leq 1$ for all $n$. Then $(f_n)_n$ is equicontinuous on every compact set of $M$, and the diagonal process ensures that there exists a subsequence $(f_{\phi(n)})_n$ which converges uniformly on every compact set of $M$ to a function $f \in \mathcal{B}_\gamma$ satisfying $N_{\infty,\gamma,\tilde{\gamma}}(f) \leq 1$. To prove (i), it suffices now to show that $\lim_n |f - f_{\phi(n)}|_{\tilde{\gamma}} = 0$. Observe that $|f - f_n|_\gamma \leq \lambda_0^{-1} N_{\infty,\gamma,\tilde{\gamma}}(f - f_n) \leq 2\lambda_0^{-1}$ (proof of Proposition 5.2). Let $\varepsilon > 0$. As $\gamma < \tilde{\gamma}$, there exists a positive constant $c$ such that, for all $n \in \mathbb{N}$ and for all $x \in M$ satisfying $d(x, x_0) > c$, we have $\frac{|f(x) - f_n(x)|}{p(x)^{\tilde{\gamma}+1}} \leq \frac{2\lambda_0^{-1} p(x)^{\gamma+1}}{p(x)^{\tilde{\gamma}+1}} \leq \varepsilon$. Besides, on the compact set $M_c = \{x : x \in M, d(x, x_0) \leq c\}$, $(f_{\phi(n)})_n$ converges uniformly to $f$; thus there exists $N \in \mathbb{N}$ such that, for all $n \geq N$ and all $x \in M_c$, we have $\frac{|f(x) - f_{\phi(n)}(x)|}{p(x)^{\tilde{\gamma}+1}} \leq \varepsilon$. Consequently, for $n \geq N$, we obtain $|f - f_{\phi(n)}|_{\tilde{\gamma}} \leq \varepsilon$.

(ii)  Now let $(f_n)_n$ be a sequence of functions in $\mathcal{B}_\gamma$ such that $N_{1,\gamma}(f_n) \leq 1$. Since $N_{1,\gamma}$ and $N_{\infty,\gamma}$ are equivalent (Proposition 5.2), the sequence $(f_n)_n$ is bounded in $(\mathcal{B}_\gamma, N_{\infty,\gamma})$ by a constant $c'$. As above, we can check that there exists a subsequence $(f_{\phi(n)})_n$ which converges pointwise to a function $f \in \mathcal{B}_\gamma$. Since $|f_n| \leq c' p^{\gamma+1}$ and $p^{\gamma+1}$ is $\nu$-integrable, the Lebesgue theorem ensures that $\lim_n \nu(|f - f_{\phi(n)}|) = 0$.  $\square$

5.3. *Quasi-compactness of $P$ on $\mathcal{B}_\gamma$.*   The following statement shows that, for $\gamma \in \,]0, \gamma_0]$, $P$ is a quasi-compact operator on $\mathcal{B}_\gamma$. This property will also follow from arguments given in Section 7, but Theorem 5.5 provides a precise description of the peripheral spectrum of $P$: 1 is a simple eigenvalue and it is the unique peripheral spectral value of $P$.

THEOREM 5.5.   *For every $\gamma \in \,]0, \gamma_0]$, $P$ is a bounded linear operator on $\mathcal{B}_\gamma$, and we have the following decomposition:*

$$\mathcal{B}_\gamma = (\mathbb{C} \cdot 1) \oplus H_\gamma,$$

*where $H_\gamma = \{f : f \in \mathcal{B}_\gamma, \; \nu(f) = 0\}$ is a closed $P$-invariant subspace of $\mathcal{B}_\gamma$ such that $r(P_{|H_\gamma}) \leq (\vartheta_0)^{1/n_0} < 1$; the real number $\vartheta_0 < 1$ has been defined in Lemma 5.1, and $r(P_{|H_\gamma})$ is the spectral radius of the restriction of $P$ to $H_\gamma$.*

PROOF.   Here it is convenient to consider $\mathcal{B}_\gamma$ equipped with the norm $N_\gamma$. We have, for all $k \geq 1$,

$$\int_G \frac{\Delta_\gamma(gx, gy)}{\Delta_\gamma(x, y)} \, d\pi^{*k}(g) = \int_G \frac{d(gx, gy)}{d(x, y)} \left( \frac{p(gx)}{p(x)} \right)^\gamma \left( \frac{p(gy)}{p(y)} \right)^\gamma d\pi^{*k}(g)$$

$$\leq \int_G c(g)\delta(g)^{2\gamma} \, d\pi^{*k}(g) = D_k(\gamma).$$

Since $\delta \leq \tilde{\delta}$, and $c$ and $\tilde{\delta}$ are submultiplicative (Lemma 4.1), hypothesis $\mathcal{M}'_{2\gamma_0 + 1} < +\infty$ and Fubini's theorem ensure that $D_k(\gamma) < +\infty$. Let $f \in \mathcal{B}_\gamma$. We have, for $x, y \in M$,

$$|P^k f(x) - P^k f(y)| \leq \int_G |f(gx) - f(gy)| \, d\pi^{*k}(g)$$

$$\leq m_\gamma(f)\Delta_\gamma(x, y) \int_G \frac{\Delta_\gamma(gx, gy)}{\Delta_\gamma(x, y)} \, d\pi^{*k}(g)$$

$$\leq m_\gamma(f)\Delta_\gamma(x, y)D_k(\gamma).$$

With $k = 1$, the foregoing proves that $Pf \in \mathcal{B}_\gamma$, and $m_\gamma(Pf) \leq D_1(\gamma) \, m_\gamma(f)$. Since $\nu(Pf) = \nu(f)$, we see that $P$ is a bounded linear operator on $(\mathcal{B}_\gamma, N_\gamma)$. As $\nu(p^{\gamma+1}) < +\infty$, the distribution $\nu$ defines a continuous linear functional on $\mathcal{B}_\gamma$; consequently, $H_\gamma = \operatorname{Ker}\nu$ is a closed subspace; it is $P$-invariant because $\nu P = \nu$.

On the other hand, with $k = n_0$, since $D_{n_0}(\gamma) \leq \vartheta_0$ (Lemma 5.1), we get $m_\gamma(P^{n_0} f) \leq \vartheta_0 m_\gamma(f)$, and by induction, $m_\gamma(P^{qn_0} f) \leq \vartheta_0^q m_\gamma(f)$ for every $q \geq 0$. In particular, if $h \in H_\gamma$, then, for every $q \geq 1$, we have $\nu(P^{qn_0} h) = \nu(h) = 0$, thus $N_\gamma(P^{qn_0} h) = m_\gamma(P^{qn_0} h) \leq \vartheta_0^q m_\gamma(h) = \vartheta_0^q N_\gamma(h)$. Thus $r(P_{|H_\gamma}) = (r(P^{n_0}_{|H_\gamma}))^{1/n_0} \leq (\vartheta_0)^{1/n_0}$.

The identity $f = \nu(f) \cdot 1 + (f - \nu(f) \cdot 1)$ leads to the stated decomposition.   □

**6. Fourier operators on $\mathcal{B}_\gamma$.** Recall that the Fourier kernels $P(t)$, $t \in \mathbb{R}$, associated to $P$ and $\xi$ are defined by

$$(P(t)f)(x) = \int_G e^{it\xi(g,x)} f(gx) \, d\pi(g),$$

and that $\xi$ is a real-valued function on $G \times M$ satisfying Condition RS.

We shall prove that, for suitable $\eta'$, $P(t)$ acts continuously on $\mathcal{B}_{\eta'}$. But, for $0 < \eta' < \eta$, it will also be convenient to see $P(t)$ as a bounded linear map from $\mathcal{B}_{\eta'}$ to $\mathcal{B}_\eta$; this is true by virtue of the following topological embedding that will be exploited repeatedly in the sequel:

> if $0 < \eta' < \eta$, we have $\mathcal{B}_{\eta'} \subset \mathcal{B}_\eta$ and for all $f \in \mathcal{B}_{\eta'}$,
>
> $N_{\infty,\eta}(f) \le N_{\infty,\eta'}(f)$.

Let $\mathcal{L}(\mathcal{B}_{\eta'}, \mathcal{B}_\eta)$ be the space of all bounded linear maps from $(\mathcal{B}_{\eta'}, N_{\infty,\eta'})$ to $(\mathcal{B}_\eta, N_{\infty,\eta})$. We denote by $\| \cdot \|_{\eta',\eta}$ the operator norm on $\mathcal{L}(\mathcal{B}_{\eta'}, \mathcal{B}_\eta)$; when $\eta' = \eta$, we merely set $\| \cdot \|_\eta = \| \cdot \|_{\eta,\eta}$.

6.1. *Preliminary remarks about the function $P(\cdot)$.* As already mentioned in Section 2, the spectral method described in Hennion and Hervé (2001) consists in applying perturbation theory to $P(t)$, so that the map $P(\cdot)$ has to be sufficiently regular.

In order to understand what are the restrictions imposed here by this property, suppose that Condition RS holds with $r > 0$, and let us study the quantity $|P(t)f - Pf|_\gamma$ for $f \in \mathcal{B}_\gamma$. Let $\varepsilon \in ]0,1]$. From the inequality $|e^{iu} - 1| \le 2|u|^\varepsilon$, Condition RS, and Lemma 4.1, we have, for all $x \in M$,

$$|P(t)f(x) - Pf(x)|$$

$$\le \int_G |e^{it\xi(g,x)} - 1| \, |f(gx)| \, d\pi(g)$$

$$\text{(I')} \qquad \le 2|t|^\varepsilon (1 + d(x,x_0))^{r\varepsilon} |f|_\gamma \, p(x)^{\gamma+1} \int_G R(g)^\varepsilon \frac{p(gx)^{\gamma+1}}{p(x)^{\gamma+1}} \, d\pi(g)$$

$$\text{(I)} \qquad \le 2C|t|^\varepsilon |f|_\gamma \, p(x)^{\gamma+1} (1 + d(x,x_0))^{r\varepsilon}$$

with $C = \int_G R(g)^\varepsilon \delta(g)^{\gamma+1} \, d\pi(g)$. Because $(1 + d(\cdot, x_0))^{r\varepsilon}$ is not bounded on $M$, this estimation does not imply that $\lim_{t \to 0} |P(t)f - Pf|_\gamma = 0$. Similar complications appear when one considers $m_\gamma(P(t)f - Pf)$.

To get around these difficulties in the special case of autoregressive processes (Section 3), Milhaud and Raugi (1989) have used a space of locally Lipschitz functions similar to $\mathcal{B}_\gamma$, which is defined by replacing $p(\cdot)^{\gamma+1}$ with $p(\cdot)^{\gamma+1} e^{\lambda d(\cdot, x_0)}$, where $\lambda$ is a positive parameter. In this case, provided the strict contraction and exponential moment conditions given in the above mentioned paper are satisfied, one can verify that the right member of $(I')$ is bounded, and more generally, that $P(\cdot)$ is a regular function from a neighborhood of $t = 0$ to $\mathcal{L}(\mathcal{B}_\gamma)$.

In this paper, we use another method which enables us to weaken the contraction and moment hypotheses considered in previous papers. This method is based on the two following facts:

1. By integrating (I) with respect to the measure $\nu$, we obtain $\nu(|P(t)f - Pf|) \le C'|t|^\varepsilon N_{1,\gamma}(f)$. This weak continuity property will be sufficient to apply to $P(t)$ a perturbation theorem of Keller and Liverani (1999).

2. Let $0 < \eta' < \eta$. For $f \in \mathcal{B}_{\eta'}$, we have

$$|P(t)f(x) - Pf(x)| \le 2C\lambda_0^{-r\varepsilon}|t|^\varepsilon |f|_{\eta'} p(x)^{\eta' + r\varepsilon + 1},$$

so that, if $\eta' + r\varepsilon \le \eta$, we get $|P(t)f - Pf|_\eta \le 2C\lambda_0^{-r\varepsilon}|t|^\varepsilon |f|_{\eta'}$. This leads us to investigate the continuity and, more generally, the existence of the Taylor expansions of $P(t)$ at $t = 0$ when $P(\cdot)$ is viewed as an $\mathcal{L}(\mathcal{B}_{\eta'}, \mathcal{B}_\eta)$-valued map [instead of an $\mathcal{L}(\mathcal{B}_{\eta'})$-valued map]; this is the aim of Sections 6.2 and 6.3. Let us mention that similar methods are used in Le Page (1989) and Hennion (1991) for other purposes.

6.2. *Taylor expansions of $P(t)$ at $t = 0$.* For $\tau > 0$ and any nonnegative measurable functions $U$, $V$ on $G$, we set

$$\mathcal{I}^\tau(U, V) = \int_G U(g)\, c(g)\, \delta(g)^{2\tau}\, d\pi(g) + \int_G V(g)\, \delta(g)^{\tau + 1}\, d\pi(g).$$

$\mathcal{I}^\tau(U, V)$ is an additive positively homogeneous function of both $U$ and $V$, and an increasing function of the variable $\tau$ because $\delta(\cdot) \ge 1$.

Observe that, for $0 < \gamma \le \gamma_0$, we have $\mathcal{I}^\gamma(1, 1) \le \mathcal{M}'_{2\gamma_0 + 1} + \mathcal{M}_{\gamma_0 + 1} < +\infty$. Let us state the three main results of this section.

PROPOSITION 6.1. *Suppose $s + 1 \le \gamma_0$, and let $\gamma$ be a real number such that $s + 1 \le \gamma \le \gamma_0$ and*

$$\mathcal{I}^\gamma(0, S) = \int_G S(g)\delta(g)^{\gamma + 1}\, d\pi(g) < +\infty.$$

*Then, for all $t \in \mathbb{R}$, $P(t) \in \mathcal{L}(\mathcal{B}_\gamma)$. Besides, there exists a constant $C$ such that we have, for all $f \in \mathcal{B}_\gamma$,*

$$|P(t)^{n_0} f|_\gamma \le \mathcal{I}^\gamma(0, 1)\, |f|_\gamma, \qquad m_\gamma(P(t)^{n_0} f) \le \vartheta_0 m_\gamma(f) + C|t|\mathcal{I}^\gamma(0, S)|f|_\gamma,$$

*where $\vartheta_0 < 1$ is the real number defined in Lemma 5.1.*

PROPOSITION 6.2. *Suppose that the following condition holds:*

$$\mathcal{U}_0(\eta', \eta)\!: 0 < \eta' \le \gamma_0,\ \eta' < \eta,\ s + 1 \le \eta,\ \mathcal{I}^{\eta'}(0, S) < +\infty.$$

*Then*

$$\lim_{|t| \to 0} \|P(t) - P\|_{\eta', \eta} = 0.$$

With the view of obtaining the Taylor expansions of $P(t)$ at $t = 0$, let us introduce, for $k \in \mathbb{N}^*$, the kernels

$$(L_k f)(x) = \int_G (i\xi(g, x))^k f(gx) \, d\pi(g).$$

PROPOSITION 6.3.   *Let $n \geq 1$. Suppose that the following condition holds:*

$$\mathcal{U}_n(\eta', \eta) : 0 < \eta' \leq \gamma_0, \ \eta' + nr < \eta, \ s + 1 + (n-1)r < \eta,$$

$$\mathcal{I}^{\eta'}(R^n, (R+S)R^{n-1}) < +\infty.$$

Then, for $k = 1, \ldots, n$, $\ L_k \in \mathcal{L}(\mathcal{B}_{\eta'}, \mathcal{B}_\eta)$ and

$$\lim_{|t| \to 0} \frac{1}{|t|^n} \left\| P(t) - P - \sum_{k=1}^n \frac{t^k}{k!} L_k \right\|_{\eta', \eta} = 0.$$

6.3.  *Proofs of Propositions* 6.1–6.3.   The main tool is Lemma 6.4, which will be stated in the next technical context.

Let $k \in \mathbb{N}^*$. Consider a complex-valued measurable function $q$ on $G^k \times M$.

Let $\alpha, \beta \in \mathbb{R}_+$, and let $A, B$ be nonnegative measurable functions on $G^k \times M$. We shall say that the inequalities (A) and (B) are satisfied if, for all $h \in G^k$ and for all $(x, y) \in M^2$ satisfying $d(x_0, y) \leq d(x_0, x)$, we have:

(A)  $|q(h, x)| \leq A(h, x) p(x)^\alpha$,
(B)  $|q(h, x) - q(h, y)| \leq B(h, x) \, d(x, y) p(x)^\beta$.

For $x \in M$, we denote by $A_x$ and $B_x$ the nonnegative functions defined on $G^k$ by $A_x(h) = A(h, x)$ and $B_x(h) = B(h, x)$.

For $h = (h_1, \ldots, h_k) \in G^k$, we set $h^\star = h_1 \cdots h_k$, and we denote by $\pi^{\otimes k}$ the product measure on $G^k$. If $x \in M$ and if $f$ is a measurable function on $M$ such that $h \mapsto q(h, x) f(h^\star x)$ is $\pi^{\otimes k}$-integrable, then we set

$$(\mathcal{K}f)(x) = \int_{G^k} q(h, x) f(h^\star x) \, d\pi^{\otimes k}(h).$$

For $\tau > 0$ and for any nonnegative measurable functions $U, V$ on $G^k$, we set

$$\mathcal{I}_k^\tau(U, V) = \int_{G^k} U(h) c(h^\star) \delta(h^\star)^{2\tau} \, d\pi^{\otimes k}(h) + \int_{G^k} V(h) \delta(h^\star)^{\tau+1} \, d\pi^{\otimes k}(h).$$

This integral only occurs in the following technical lemma; notice that it equals $\mathcal{I}^\tau(U, V)$ when $k = 1$.

Lemma 6.4.   *Let $0 < \eta' \leq \eta$. Suppose that, for all $x \in M$, we have $\mathcal{I}_k^{\eta'}(A_x, A_x + B_x) < +\infty$. Then, for $f \in \mathcal{B}_{\eta'}$ and $x \in M$, $\mathcal{K}f(x)$ is defined; moreover, for $x, y \in M$ such that $x \neq y$ and $d(y, x_0) \leq d(x, x_0)$, we have the inequalities*

$$\frac{|\mathcal{K}f(x)|}{p(x)^{\eta+1}} \leq \frac{\mathcal{I}_k^{\eta'}(0, A_x)}{p(x)^{\eta-\eta'-\alpha}} \, |f|_{\eta'},$$

$$\frac{|\mathcal{K}f(x) - \mathcal{K}f(y)|}{\Delta_\eta(x,y)} \leq \frac{\mathcal{I}_k^{\eta'}(A_x, 0)}{p(x)^{\eta-\eta'-\alpha}} \, m_{\eta'}(f) \; + \; \frac{\mathcal{I}_k^{\eta'}(0, B_x)}{p(x)^{\eta-\beta-1}} \, |f|_{\eta'}.$$

To apply this lemma, it will be worth noticing that, for $\eta > 0$ and for any function $f$ on $M$, we have, owing to symmetry,

$$m_\eta(f) = \sup\left\{ \frac{|f(x) - f(y)|}{\Delta_\eta(x,y)}, x, y \in M, x \neq y, d(y, x_0) \leq d(x, x_0) \right\}.$$

Proof of Lemma 6.4.   We shall use the inequalities $\sup_{x \in M} p(gx)/p(x) \leq \delta(g)$ (Lemma 4.1) and $|f(g \cdot)| \leq |f|_{\eta'} p(g \cdot)^{\eta'+1}$. Let $f$, $x$ and $y$ be as in the statement. We have

$$\int_{G^k} |q(h,x) f(h^\star x)| \, d\pi^{\otimes k}(h) \leq p(x)^\alpha |f|_{\eta'} \int_{G^k} A(h,x) p(h^\star x)^{\eta'+1} \, d\pi^{\otimes k}(h)$$

$$\leq p(x)^{\alpha+\eta'+1} |f|_{\eta'} \mathcal{I}_k^{\eta'}(0, A_x).$$

It follows that $\mathcal{K}f(x)$ is defined and verifies the first stated inequality.

To prove the second one, let us write

$$|\mathcal{K}f(x) - \mathcal{K}f(y)| \leq A_1(x,y) + A_2(x,y)$$

with

$$A_1(x,y) = \int_{G^k} |q(h,x)| \, |f(h^\star x) - f(h^\star y)| \, d\pi^{\otimes k}(h),$$

$$A_2(x,y) = \int_{G^k} |f(h^\star y)| \, |q(h,x) - q(h,y)| \, d\pi^{\otimes k}(h).$$

Then

$$\frac{A_1(x,y)}{\Delta_\eta(x,y)} \leq m_{\eta'}(f) p(x)^\alpha$$

$$\times \int_{G^k} A(h,x) \frac{d(h^\star x, h^\star y) p(h^\star x)^{\eta'} p(h^\star y)^{\eta'}}{d(x,y) p(x)^\eta p(y)^\eta} \, d\pi^{\otimes k}(h)$$

$$\leq m_{\eta'}(f) \left( \frac{p(x)^\alpha}{p(x)^{\eta-\eta'}} \right) \left( \frac{1}{p(y)^{\eta-\eta'}} \right)$$

$$\times \int_{G^k} A(h,x) c(h^\star) \delta(h^\star)^{2\eta'} \, d\pi^{\otimes k}(h)$$

(M1) $$\leq \frac{\mathcal{I}_k^{\eta'}(A_x, 0)}{p(x)^{\eta - \eta' - \alpha}} m_{\eta'}(f) \qquad [\text{because } p(y)^{\eta - \eta'} \geq 1].$$

Consider now the quantity $A_2(x, y)$. By using the inequality $d(y, x_0) \leq d(x, x_0)$, we obtain

$$A_2(x, y) \leq |f|_{\eta'} d(x, y) p(x)^{\beta} \int_{G^k} B(h, x) p(h^{\star} y)^{\eta' + 1} d\pi^{\otimes k}(h)$$

$$\leq |f|_{\eta'} d(x, y) p(x)^{\beta} p(y)^{\eta' + 1} \mathcal{I}_k^{\eta'}(0, B_x)$$

(M'2) $$\leq |f|_{\eta'} d(x, y) p(x)^{\beta + 1} p(y)^{\eta'} \mathcal{I}_k^{\eta'}(0, B_x),$$

and from $p(y)^{\eta'} \leq p(y)^{\eta}$, we get

(M2) $$\frac{A_2(x, y)}{\Delta_{\eta}(x, y)} \leq \frac{\mathcal{I}_k^{\eta'}(0, B_x)}{p(x)^{\eta - \beta - 1}} |f|_{\eta'}.$$

We conclude by combining (M1) and (M2). $\quad \square$

We shall also need the next bounds.

LEMMA 6.5. *For $n \in \mathbb{N}$ and $x \in \mathbb{R}$, we set $\phi_n(x) = e^{ix} - \sum_{k=0}^{n} \frac{(ix)^k}{k!}$. For all $x, y \in \mathbb{R}$, we have:*

  (i) $|\phi_n(x)| \leq 2|x|^n \min\{1, |x|\}$,
  (ii) $|e^{iy} - e^{ix}| \leq |y - x|$,
  (iii) *for $n \geq 1$, $|\phi_n(y) - \phi_n(x)| \leq 2|y - x|(|x|^{n-1} \min\{1, |x|\} + |y|^{n-1} \times \min\{1, |y|\})$.*

PROOF. The assertion (ii) is clear, and it implies that $|\phi_0(x)| \leq \min\{2, |x|\} \leq 2 \min\{1, |x|\}$.

Let $n \geq 1$. The Taylor formula to the orders $n$ and $n - 1$ with integral remainder shows that

$$|\phi_n(x)| \leq \frac{|x|^{n+1}}{(n+1)!} \quad \text{and} \quad |\phi_n(x)| = \left| \phi_{n-1}(x) - \frac{(ix)^n}{n!} \right| \leq 2 \frac{|x|^n}{n!}.$$

Hence

$$|\phi_n(x)| \leq \min \left\{ \frac{2|x|^n}{n!}, \frac{|x|^{n+1}}{(n+1)!} \right\} \leq 2|x|^n \min\{1, |x|\}.$$

Since $\phi_n'(x) = i\phi_{n-1}(x)$ for $n \geq 1$, we have $|\phi_n(x) - \phi_n(y)| \leq |x - y| \times \sup\{|\phi_{n-1}(t)| : t \in [x, y]\}$. This inequality and point (i) prove assertion (iii). $\quad \square$

Now let us prove Propositions 6.1–6.3.

PROOF OF PROPOSITION 6.1.   Let $k \in \mathbb{N}^*$. By induction, we easily prove that

$$(P(t)^k f)(x) = \int_{G^k} e^{it\xi_k(h,x)} f(h^\star x) \, d\pi^{\otimes k}(h),$$

with $\xi_k(h,x) = \xi(h_k,x) + \xi(h_{k-1},h_k x) + \cdots + \xi(h_1,h_2 \cdots h_k x)$, for all $h = (h_1, \ldots, h_k) \in G^k$. Therefore $\mathcal{K} = P(t)^k$ is associated to the kernel $q(h,x) = e^{it\xi_k(h,x)}$. We have $|q(h,x)| = 1$; Condition RS and Lemma 4.1 give, for $g_1, g_2 \in G$,

$$|\xi(g_1,g_2x) - \xi(g_1,g_2y)| \leq S(g_1) \, d(g_2x,g_2y)(1 + d(g_2x,x_0) + d(g_2y,x_0))^s$$
$$\leq \lambda_0^{-s} \, d(x,y) S(g_1) c(g_2)(p(g_2x) + p(g_2y))^s$$
$$\leq \lambda_0^{-s} \, d(x,y) S(g_1) c(g_2) \delta(g_2)^s (p(x) + p(y))^s.$$

Hence, if $d(x_0,y) \leq d(x_0,x)$, we get

$$|\xi(g_1,g_2x) - \xi(g_1,g_2y)| \leq 2^s \lambda_0^{-s} \, d(x,y) \, S(g_1) \, c(g_2) \delta(g_2)^s p(x)^s.$$

Finally, by using Lemma 6.5(ii) and the facts that $\delta(\cdot) \leq \tilde{\delta}(\cdot)$ and that $c(\cdot)$, $\tilde{\delta}(\cdot)$ are submultiplicative (Lemma 4.1), we obtain that $q(h,x)$ verifies the inequalities (A) and (B) with

$$A(h,x) = 1, \qquad \alpha = 0, \qquad B(h,x) = 2^s \lambda_0^{-s} |t| B_k(h), \qquad \beta = s,$$

where $B_k(h) = \sum_{i=1}^k S(h_i) c(h_{i+1}) \cdots c(h_k) \tilde{\delta}(h_{i+1})^s \cdots \tilde{\delta}(h_k)^s$. We have

$$\mathcal{I}_k^\eta(A, A + B) = \mathcal{I}_k^\eta(1, 1 + 2^s \lambda_0^{-s} |t| B_k) \leq \mathcal{I}_k^{\gamma_0}(1,1) + 2^s \lambda_0^{-s} |t| \mathcal{I}_k^{\gamma_0}(0, B_k).$$

Since $c\delta^{2\gamma_0} \leq c\tilde{\delta}^{2\gamma_0}$, $\delta^{\gamma_0+1} \leq \tilde{\delta}^{\gamma_0+1}$, and the functions $c$, $\tilde{\delta}$ are submultiplicative, hypotheses $\mathcal{M}_{\gamma_0+1} < +\infty$, $\mathcal{M}'_{2\gamma_0+1} < +\infty$, and Fubini's theorem imply that $\mathcal{I}_k^{\gamma_0}(1,1) < +\infty$. Besides, we have

$$\mathcal{I}_k^\gamma(0, B_k) \leq \sum_{i=1}^k \int_{G^k} S(h_i) c(h_{i+1}) \cdots c(h_k) \tilde{\delta}(h_{i+1})^s$$
$$\cdots \tilde{\delta}(h_k)^s \, \tilde{\delta}(h_1)^{\gamma+1} \cdots \tilde{\delta}(h_k)^{\gamma+1} \, d\pi^{\otimes k}(h).$$

We have $\tilde{\delta} \leq \frac{2}{\lambda_0} \delta$, and $\int_G S(g) \delta(g)^{\gamma+1} \, d\pi(g) < +\infty$, thus $\int_G S(g) \times \tilde{\delta}(g)^{\gamma+1} \, d\pi(g) < +\infty$. Moreover, we have $c(g)\tilde{\delta}(g)^{\gamma+1+s} \leq c(g)\tilde{\delta}(g)^{2\gamma_0}$. It follows from hypothesis $\mathcal{M}'_{2\gamma_0+1} < +\infty$ and Fubini's theorem that $\mathcal{I}_k^\gamma(0, B_k) < +\infty$.

Now let us apply Lemma 6.4 with $\eta' = \eta = \gamma \leq \gamma_0$.

For $k = 1$, we get, for all $f \in \mathcal{B}_\gamma$,

$$|P(t)f|_\gamma \leq \mathcal{I}^\gamma(0,1) \, |f|_\gamma.$$

On the other hand, since $\gamma \geq s + 1$, we have $p(x)^{\gamma-s-1} \geq 1$; hence, since $B_1 = S$,

$$m_\gamma(P(t)) \leq \mathcal{I}^\gamma(1,0) \, m_\gamma(f) + 2^s \lambda_0^{-s} |t| \mathcal{I}^\gamma(0,S) |f|_\gamma.$$

This proves that $P(t) \in \mathcal{L}(\mathcal{B}_\gamma)$.

For $k = n_0$, the first inequality is still valid for $P(t)^{n_0}$, while the second one becomes

$$m_\gamma(P(t)^{n_0} f) \leq \mathcal{I}_{n_0}^\gamma(1, 0) m_\gamma(f) + 2^s \lambda_0^{-s} |t| \mathcal{I}_{n_0}^\gamma(0, B_{n_0}) |f|_\gamma,$$

with $\mathcal{I}_{n_0}^\gamma(1, 0) \leq \mathcal{I}_{n_0}^{\gamma_0}(1, 0) = \int_G c(h) \delta(h)^{2\gamma_0} \, d\pi^{*n_0}(h) = \vartheta_0.$   $\square$

To establish Propositions 6.2 and 6.3, we shall employ the notation

$$\tau(t, g, x) = \min\{1, |t| R(g)(1 + d(x_0, x))^r\}.$$

LEMMA 6.6.   *Let $\eta > 0$ and let $U$, $V$ be nonnegative measurable functions on $G$ such that $\mathcal{I}^\eta(U, V) < +\infty$. Then, for all $\varepsilon > 0$,*

$$\lim_{|t| \to 0} \left( \sup_{x \in M} \frac{\mathcal{I}^\eta(U(\cdot)\tau(t, \cdot, x), V(\cdot)\tau(t, \cdot, x))}{(1 + d(x, x_0))^\varepsilon} \right) = 0.$$

PROOF.   Let $\rho > 0$. We have $\tau \leq 1$ and, for $1 + d(x_0, x) \leq \rho$, we can write $\tau(t, g, x) \leq \min\{1, |t| R(g)\rho^r\} = \tau_\rho(t, g)$. Therefore, comparing $\rho$ with $1 + d(x_0, x)$, we obtain, for all $x \in M$,

$$\frac{\mathcal{I}^\eta(U(\cdot)\tau(t, \cdot, x), V(\cdot)\tau(t, \cdot, x))}{(1 + d(x, x_0))^\varepsilon}$$
$$\leq \rho^{-\varepsilon} \mathcal{I}^\eta(U, V) + \mathcal{I}^\eta(U(\cdot)\tau_\rho(t, \cdot), V(\cdot)\tau_\rho(t, \cdot)).$$

Since $\lim_{|t| \to 0} \tau_\rho(t, g) = 0$ and $\tau_\rho \leq 1$, the dominated convergence theorem implies that

$$\limsup_{|t| \to 0} \left( \sup_{x \in M} \frac{\mathcal{I}^\eta(U(\cdot)\tau(t, \cdot, x), V(\cdot)\tau(t, \cdot, x))}{(1 + d(x, x_0))^\varepsilon} \right) \leq \rho^{-\varepsilon} \mathcal{I}^\eta(U, V).$$

Since $\rho$ is arbitrary, this provides the desired statement.   $\square$

PROOF OF PROPOSITION 6.2.   Let us consider the kernel $q(g, x) = e^{it\xi(g, x)} - 1$, $x \in M$, $g \in G$, which defines the operator $\mathcal{K}_t = P(t) - P$. By Lemma 6.5(i) with $n = 0$, and then (ii), we have, for $x \in M$ and $g \in G$,

$$|q(g, x)| = |e^{it\xi(g, x)} - 1| \leq 2\min\{1, |t| |\xi(g, x)|\} \leq 2\tau(t, g, x),$$

and, if $d(y, x_0) \leq d(x, x_0)$,

$$|q(g, x) - q(g, y)| \leq |t| |\xi(g, x) - \xi(g, y)| \leq 2^s |t| S(g) \, d(x, y)(1 + d(x_0, x))^s.$$

Lemma 6.4 applied with $k = 1$, and

$$A(g, x) = 2\tau(t, g, x), \qquad \alpha = 0, \qquad B(g, x) = 2^s \lambda_0^{-s} |t| S(g), \qquad \beta = s,$$

yields

$$\frac{|\mathcal{K}_t f(x)|}{p(x)^{\eta+1}} \le 2\lambda_0^{\eta'-\eta}\, \frac{\mathcal{I}^{\eta'}(0,\tau(t,\cdot,x))}{(1+d(x,x_0))^{\eta-\eta'}}\,|f|_{\eta'},$$

$$\frac{|\mathcal{K}_t f(x) - \mathcal{K}_t f(y)|}{\Delta_\eta(x,y)} \le 2\lambda_0^{\eta'-\eta}\, \frac{\mathcal{I}^{\eta'}(\tau(t,\cdot,x),0)}{(1+d(x,x_0))^{\eta-\eta'}}\,m_{\eta'}(f)$$

$$+\, 2^s \lambda_0^{1-\eta}|t|\, \frac{\mathcal{I}^{\eta'}(0,S(\cdot))}{(1+d(x,x_0))^{\eta-s-1}}\,|f|_{\eta'}.$$

Since $\eta - \eta' > 0$, $\eta - s - 1 > 0$, $\mathcal{I}^{\eta'}(0,1) \le \mathcal{M}_{\gamma_0+1} < +\infty$, and $\mathcal{I}^{\eta'}(1,0) \le \mathcal{M}'_{2\gamma_0+1} < +\infty$, we conclude by using Lemma 6.6 with $(U,V) = (0,1)$ and $(U,V) = (1,0)$. $\square$

PROOF OF PROPOSITION 6.3.    Let us consider the kernel $q(g,x) = e^{it\xi(g,x)} - \sum_{k=0}^{n} \frac{(it\xi(g,x))^k}{k!}$, and set $\mathcal{K}_t = P(t) - P - \sum_{k=1}^{n} \frac{t^k}{k!} L_k$.

Assertion (i) of Lemma 6.5 implies that we have, for $x \in M$ and $g \in G$,

$$|q(g,x)| \le 2|t|^n |\xi(g,x)|^n \min\{1,|t|\,|\xi(g,x)|\}$$
$$\le 2|t|^n R(g)^n (1+d(x_0,x))^{nr}\tau(t,g,x),$$

while assertion (iii) shows that, for $d(x_0,y) \le d(x_0,x)$,

$$|q(g,x) - q(g,y)|$$
$$\le 2|t|^n |\xi(g,x) - \xi(g,y)|(|\xi(g,x)|^{n-1}\min\{1,\,|t|\,|\xi(g,x)|\}$$
$$+\, |\xi(g,y)|^{n-1}\min\{1,\,|t|\,|\xi(g,y)|\})$$
$$\le 2^{s+2}|t|^n S(g)\, d(x,y) R(g)^{n-1}(1+d(x_0,x))^{s+(n-1)r}\tau(t,g,x).$$

Therefore inequalities (A) and (B) hold with $k=1$, and

$$A(g,x) = 2\lambda_0^{-nr}|t|^n R^n(g)\tau(t,g,x), \qquad\qquad \alpha = nr,$$
$$B(g,x) = 2^{s+2}\lambda_0^{-\beta}|t|^n S(g)R^{n-1}(g)\tau(t,g,x), \qquad \beta = s+(n-1)r.$$

From Lemma 6.4 with $k=1$, it follows that

$$\frac{|\mathcal{K}_t f(x)|}{p(x)^{\eta+1}} \le 2\lambda_0^{\eta'-\eta}|t|^n \frac{\mathcal{I}^{\eta'}(0,R^n(\cdot)\tau(t,\cdot,x))}{(1+d(x,x_0))^{\eta-\eta'-nr}}\,|f|_{\eta'},$$

$$\frac{|\mathcal{K}_t f(x) - \mathcal{K}_t f(y)|}{\Delta_\eta(x,y)} \le 2\lambda_0^{\eta'-\eta}|t|^n \frac{\mathcal{I}^{\eta'}(R^n(\cdot)\tau(t,\cdot,x),0)}{(1+d(x,x_0))^{\eta-\eta'-nr}}\,m_{\eta'}(f)$$

$$+\, 2^{s+2}\lambda_0^{1-\eta}|t|^n \frac{\mathcal{I}^{\eta'}(0,S(\cdot)R^{n-1}(\cdot)\tau(t,\cdot,x))}{(1+d(x,x_0))^{\eta-\beta-1}}\,|f|_{\eta'}.$$

Since $\eta - \eta' - rn > 0$, $\eta - \beta - 1 = \eta - s - (n-1)r - 1 > 0$, and $\mathcal{I}^{\eta'}(R^n, (R + S)R^{n-1}) < +\infty$, the previous inequalities imply that $\mathcal{K}_t \in \mathcal{L}(\mathcal{B}_{\eta'}, \mathcal{B}_\eta)$; then, by using Lemma 6.6, we get

$$\lim_{|t| \to 0} \frac{1}{|t|^n} \|\mathcal{K}_t\|_{\eta', \eta} = 0.$$

Finally, it remains to prove that, for $k = 1, \ldots, n$, $L_k \in \mathcal{L}(\mathcal{B}_{\eta'}, \mathcal{B}_\eta)$. This derives from the following: on the one hand, $P, P(t) \in \mathcal{L}(\mathcal{B}_{\eta'}, \mathcal{B}_\eta)$ (Proposition 6.2), and on the other hand, by the above, we have $P(t) - P - \sum_{k=1}^{n'} \frac{t^k}{k!} L_k \in \mathcal{L}(\mathcal{B}_{\eta'}, \mathcal{B}_\eta)$ for $n' = 0, \ldots, n$. $\square$

To end this section, we give an additional statement which completes Proposition 6.1 and will be helpful in the proof of Proposition 7.4.

PROPOSITION 6.7. *Assume $s + 1 < \gamma_0$, and let $\eta$ and $\tilde{\eta}$ be real numbers such that $s + 1 + (\tilde{\eta} - \eta) \leq \eta < \tilde{\eta} < \gamma_0$ and*

$$\mathcal{I}^{\tilde{\eta}}(0, S) < +\infty.$$

*Then there exists a constant $C$ such that we have, for all $t \in \mathbb{R}$ and $f \in \mathcal{B}_\eta$,*

$$m_\eta(P(t)^{n_0} f) \leq \vartheta_0 \, m_\eta(f) + C|t| \, |f|_{\tilde{\eta}}.$$

PROOF. First, we establish the following with the notation of Lemma 6.4.

LEMMA 6.8. *Suppose that inequalities (A) and (B) hold. Let $0 < \eta < \tilde{\eta}$. If $\alpha = 0$, $\beta + 1 + \tilde{\eta} < 2\eta$, and if, for all $x \in M$, $\mathcal{I}_k^\eta(A_x, 0) + \mathcal{I}_k^\eta(0, B_x) < +\infty$, then we have, for all $f \in \mathcal{B}_\eta$,*

$$m_\eta(\mathcal{K}f) \leq \mathcal{I}_k^\eta(A_x, 0) m_\eta(f) + \mathcal{I}_k^{\tilde{\eta}}(0, B_x) \, |f|_{\tilde{\eta}}.$$

PROOF. Let us write, as in the proof of Lemma 6.4, $|\mathcal{K}f(x) - \mathcal{K}f(y)| \leq A_1(x, y) + A_2(x, y)$, and let us return to inequalities (M1) and (M'2).

With $\eta' = \eta$ and $\alpha = 0$, (M1) gives $\frac{A_1(x,y)}{\Delta_\eta(x,y)} \leq \mathcal{I}_k^\eta(A_x, 0) m_\eta(f)$.

Inequality (M'2) holds for any $\eta' > 0$; in particular, it is satisfied with $\eta' = \tilde{\eta}$. Besides, if $d(y, x_0) \leq d(x, x_0)$, we have $p(y)^{\tilde{\eta}} = p(y)^{\tilde{\eta} - \eta} p(y)^\eta \leq p(x)^{\tilde{\eta} - \eta} p(y)^\eta$. Hence $A_2(x, y) \leq |f|_{\tilde{\eta}} d(x, y) p(x)^{\beta + 1 + \tilde{\eta} - \eta} p(y)^\eta \mathcal{I}_k^{\tilde{\eta}}(0, B_x)$. Since $\beta + 1 + \tilde{\eta} - \eta < \eta$, we obtain

$$\frac{A_2(x, y)}{\Delta_\eta(x, y)} \leq \mathcal{I}_k^{\tilde{\eta}}(0, B_x) \, |f|_{\tilde{\eta}}.$$

We conclude by combining the two previous bounds. $\square$

Let us now prove the proposition. Consider the kernel $q(g,x) = e^{it\xi_{n_0}(h,x)}$, $h \in G^{n_0}$, $x \in M$, defining $P(t)^{n_0}$ (see proof of Proposition 6.1); it verifies inequalities (A) and (B) with $k = n_0$, and $\alpha, \beta, A, B$ given in the proof of Proposition 6.1. Lemma 6.8 applies to this kernel because $\beta + 1 + \tilde{\eta} = s + 1 + \tilde{\eta} \le 2\eta$, $\mathcal{I}_{n_0}^\eta(A_x, 0) = \mathcal{I}_{n_0}^{\gamma_0}(1,0) \le \vartheta_0$, and $\mathcal{I}_{n_0}^{\tilde{\eta}}(0, B_x) < +\infty$; this last point can be shown by using hypothesis $\mathcal{I}^{\tilde{\eta}}(0,S) < +\infty$ and a method similar to that employed in the proof of Proposition 6.1. This proves the proposition. □

## 7. The spectrum of $P(t)$ acting on $\mathcal{B}_\gamma$.

We use the standard notation $\sigma(T)$ and $r(T)$ to name the spectrum and the spectral radius of an operator $T$ [see Dunford and Schwartz (1958)]. We denote by $\mathcal{B}_\gamma'$ the topological dual space of $\mathcal{B}_\gamma$, and by $\langle \cdot, \cdot \rangle$ the canonical bilinear functional on $\mathcal{B}_\gamma' \times \mathcal{B}_\gamma$.

For $\gamma \le \gamma_0$, the $P$-invariant probability distribution $\nu$ defines an element of $\mathcal{B}_\gamma'$, and Theorem 5.5 shows that $P \in \mathcal{L}(\mathcal{B}_\gamma)$, that

$$\sigma(P) \subset \{1\} \cup \{z : z \in \mathbb{C}, |z| \le \kappa_0\} \qquad \text{with } \kappa_0 = \vartheta_0^{1/n_0} < 1,$$

and that there exists $N_{(\gamma)} \in \mathcal{L}(\mathcal{B}_\gamma)$, with spectral radius $r(N_{(\gamma)}) \le \kappa_0 < 1$, such that, for $n \ge 1$ and $f \in \mathcal{B}_\gamma$,

$$P^n f = \langle \nu, f \rangle 1 + N_{(\gamma)}^n f.$$

The following statement, which is obtained by applying to $P(\cdot)$ a perturbation theorem of Keller and Liverani (1999), asserts first that, for small $|t|$, the spectrum of $P(t)$ is close to that of $P$; second, that a spectral decomposition of the preceding type is still valid for $P(t)$; and third, that the resolvents are uniformly bounded in $t$ for $z$ ranging outside a neighborhood of the spectrum of $P$.

We shall use the following notation. Let $\kappa_0'$ and $\kappa_0''$ be real numbers such that $0 < \kappa_0 < \kappa_0' < \kappa_0'' < 1$. Let $\mathcal{D}_0$ and $\mathcal{D}_1$ be the open discs of the complex plane defined by

$$\mathcal{D}_0 = \{z : z \in \mathbb{C}, |z| < \kappa_0'\}, \qquad \mathcal{D}_1 = \{z : z \in \mathbb{C}, |z - 1| < 1 - \kappa_0''\}.$$

We denote by $\Gamma_0$ and $\Gamma_1$ the oriented circles defined, respectively, as the boundaries of $\mathcal{D}_0$ and of $\mathcal{D}_1^c$. We set

$$\mathcal{R} = \mathbb{C} \setminus (\mathcal{D}_0 \cup \mathcal{D}_1) = \{z : z \in \mathbb{C}, |z| \ge \kappa_0', \ |z - 1| \ge 1 - \kappa_0''\}.$$

PROPOSITION 7.1. *Assume that* $s + 1 \le \gamma_0$. *Let* $\gamma$ *be such that* $s + 1 \le \gamma \le \gamma_0$ *and*

$$\mathcal{I}^\gamma(0,S) = \int_G S(g)\delta(g)^{\gamma+1} \, d\pi(g) < +\infty.$$

*Then, for all* $t \in \mathbb{R}$, $P(t) \in \mathcal{L}(\mathcal{B}_\gamma)$. *Moreover, there exists an open interval* $I_\gamma$ *containing* $t = 0$ *such that we have the following spectral properties, for* $t \in I_\gamma$, *and for* $P(t)$ *acting on* $\mathcal{B}_\gamma$:

(a) $\sigma(P(t)) \subset \mathcal{D}_0 \cup \mathcal{D}_1$, and there exists $\lambda_{(\gamma)}(t) \in \mathbb{C}$ such that $\sigma(P(t)) \cap \mathcal{D}_1 = \{\lambda_{(\gamma)}(t)\}$,

(b) there exists a unique function $v_{(\gamma)}(t)$, belonging to $\mathcal{B}_\gamma$, such that we have $\langle \nu, v_{(\gamma)}(t) \rangle = 1$ and $P(t)v_{(\gamma)}(t) = \lambda_{(\gamma)}(t)v_{(\gamma)}(t)$,

(c) we have $M_\gamma = \sup\{\|(z - P(t))^{-1}\|_\gamma, t \in I_\gamma, z \in \mathcal{R}\} < +\infty$,

(d) there exist $\phi_{(\gamma)}(t) \in \mathcal{B}'_\gamma$ and $N_{(\gamma)}(t) \in \mathcal{L}(\mathcal{B}_\gamma)$ such that

$$\forall f \in \mathcal{B}_\gamma, \ \forall n \in \mathbb{N}^*, \qquad P(t)^n f = \lambda_{(\gamma)}(t)^n \langle \phi_{(\gamma)}(t), f \rangle v_{(\gamma)}(t) + N_{(\gamma)}(t)^n f,$$

with $\|N_{(\gamma)}(t)^n\|_\gamma \leq \frac{M_\gamma}{2\pi}(\kappa'_0)^n$.

Notice that, for $t = 0$, we have $\lambda_{(\gamma)}(0) = 1$, $v_{(\gamma)}(0) = 1$, $\phi_{(\gamma)}(0) = \nu$ and $N_{(\gamma)}(0) = N_{(\gamma)}$. From the inclusion $\mathcal{B}_{\gamma'} \subset \mathcal{B}_\gamma$, for $0 < \gamma' < \gamma$, and from Proposition 7.1, we deduce the following corollary.

COROLLARY 7.2. Under the conditions of Proposition 7.1, if $s + 1 \leq \gamma' < \gamma \leq \gamma_0$, then, for all $t \in I_{\gamma'} \cap I_\gamma$, we have

$$\lambda_{(\gamma)}(t) = \lambda_{(\gamma')}(t), \qquad v_{(\gamma)}(t) = v_{(\gamma')}(t),$$

$$\phi_{(\gamma)}(t)_{|\mathcal{B}_{\gamma'}} = \phi_{(\gamma')}(t), \qquad N_{(\gamma)}(t)_{|\mathcal{B}_{\gamma'}} = N_{(\gamma')}(t).$$

NOTATION. In accordance with this corollary, when Proposition 7.1 applies to $P(t)$ acting on $\mathcal{B}_\gamma$, we set

$$\lambda(t) = \lambda_{(\gamma)}(t), \qquad v(t) = v_{(\gamma)}(t), \qquad \phi(t) = \phi_{(\gamma)}(t), \qquad N(t) = N_{(\gamma)}(t).$$

It will follow from the proof of Proposition 7.1 that we have the following:

COROLLARY 7.2′. Under the conditions of Proposition 7.1, for $s + 1 \leq \gamma \leq \gamma_0$ and for $t \in I_\gamma$, the elements $N(t), v(t), \phi(t)$ are given by the following formulae in which integration is considered in the space $\mathcal{L}(\mathcal{B}_\gamma)$:

$$N(t) = \frac{1}{2i\pi} \int_{\Gamma_0} (z - P(t))^{-1} \, dz,$$

$$v(t) = \frac{1}{\nu(\Pi(t)1)}\Pi(t)1,$$

$$\phi(t) = \Pi(t)^* \nu,$$

where

$$\Pi(t) = \frac{1}{2i\pi} \int_{\Gamma_1} (z - P(t))^{-1} \, dz.$$

Moreover, we have

$$N(t)^n = \frac{1}{2i\pi} \int_{\Gamma_0} z^n (z - P(t))^{-1} \, dz \quad and \quad \|N(t)^n\|_\gamma \leq \frac{M_\gamma}{2\pi}(\kappa'_0)^n.$$

PROOF OF PROPOSITION 7.1. The hypotheses are those of Proposition 6.1. Consequently, for all $t \in \mathbb{R}$, $P(t) \in \mathcal{L}(\mathcal{B}_\gamma)$.

To establish the assertions (a)–(d), we shall use the results of Keller and Liverani (1999). Let us specify the context of this paper: the space (here $\mathcal{B}_\gamma$) on which the collection of operators [here $P(t)$, $t \in \mathbb{R}$] acts, is endowed with a norm (here $N_{1,\gamma}$, Section 5.2) with respect to which the space is complete, and with an auxiliary norm which is dominated by the preceding one. An easy adaption shows that the results of Keller and Liverani (1999) are still valid with an auxiliary seminorm [here $\nu(|\cdot|)$]. The lemma below proves that the required hypotheses are fulfilled.

LEMMA 7.3. *Under the hypotheses of Proposition 7.1:*

   (i) *for $t \in \mathbb{R}$, $n \in \mathbb{N}^*$ and $f \in \mathcal{B}_\gamma$, we have $\nu(|P(t)^n f|) \le \nu(|f|)$,*

   (ii) *there exist $J \in \mathbb{R}_+$ and an open interval $I_\gamma$ containing $t = 0$ such that, for $t \in I_\gamma$, we have*

$$\forall f \in \mathcal{B}_\gamma, \qquad N_{1,\gamma}(P(t)^{n_0} f) \le (\kappa'_0)^{n_0} N_{1,\gamma}(f) + J\nu(|f|),$$

   (iii) *for all $t \in I_\gamma$, the essential spectral radius of $P(t)$ is $\le \kappa'_0$,*

   (iv) *there exists a positive continuous function $\varphi$, vanishing at $t = 0$, such that we have, for all $f \in \mathcal{B}_\gamma$, $\nu(|P(t)f - Pf|) \le \varphi(t)N_{1,\gamma}(f)$.*

We refer to Hennion and Hervé [(2001), Chapter XIV] for the notion of essential spectral radius of an operator. The property (iv) above means that, in a weak sense, for small $|t|$, $P(t)$ is a perturbation of $P$.

PROOF OF LEMMA 7.3. (i) As $P$ is nonnegative, we get $|P(t)^n f| \le P^n|f|$; hence the inequality of point (i), since $\nu$ is $P$-invariant.

   (ii) From Proposition 6.1, we have, for all $f \in \mathcal{B}_\gamma$,

$$m_\gamma(P(t)^{n_0} f) \le \kappa_0^{n_0} m_\gamma(f) + C|t|\mathcal{I}^\gamma(0, S)|f|_\gamma.$$

As a consequence of the equivalence of the norms $N_{\infty,\gamma}$ and $N_{1,\gamma}$, we get a constant $K'$, such that, for all $f \in \mathcal{B}_\gamma$, we have

$$m_\gamma(P(t)^{n_0} f) \le \kappa_0^{n_0} m_\gamma(f) + K'|t|N_{1,\gamma}(f) = (\kappa_0^{n_0} + K'|t|)m_\gamma(f) + K'|t|\nu(|f|).$$

so that, for $|t| \le \frac{\kappa'^{n_0}_0 - \kappa_0^{n_0}}{K'}$, we obtain

$$m_\gamma(P(t)^{n_0} f) \le (\kappa'_0)^{n_0} m_\gamma(f) + (\kappa'^{n_0}_0 - \kappa_0^{n_0})\nu(|f|),$$

Using point (i), we get $N_{1,\gamma}(P(t)^{n_0} f) \le (\kappa'_0)^{n_0} N_{1,\gamma}(f) + J\nu(|f|)$ with $J = \kappa'^{n_0}_0 - \kappa_0^{n_0} + 1$.

(iii) Recall that the essential spectral radius of an operator is smaller than its spectral radius; consequently, point (iii) is clear when $r(P(t)) \leq \kappa_0'$.

Assume that $r(P(t)) > \kappa_0'$. Then, from point (i) and the Doeblin–Fortet inequality established in point (ii), and from the fact that the canonical embedding of $(\mathcal{B}_\gamma, N_{1,\gamma})$ into $(\mathcal{B}_\gamma, \nu(|\cdot|))$ is compact (Lemma 5.4), we deduce by means of the Ionescu-Tulcea and Marinescu theorem or more precisely of Corollary 1 in Hennion (1993) that, for $|t| \leq \frac{\kappa_0^{m_0} - \kappa_0^{n_0}}{K'}$, $P(t)$ is quasi-compact, and that its essential spectral radius is $\leq \kappa_0'$.

(iv) Using the inequality $|f(gx)| \leq |f|_\gamma p(gx)^{\gamma+1} \leq |f|_\gamma \tilde{\delta}(g)^{\gamma+1} p(x)^{\gamma+1}$ (Lemma 4.1), we obtain

$$\nu(|P(t)f - Pf|) \leq \int_G \int_M |e^{it\xi(g,x)} - 1| \, |f(gx)| \, d\pi(g) \, d\nu(x) \leq |f|_\gamma \varepsilon(t),$$

with

$$\varepsilon(t) = \int_G \int_M |e^{it\xi(g,x)} - 1| \tilde{\delta}(g)^{\gamma+1} p(x)^{\gamma+1} \, d\pi(g) \, d\nu(x).$$

Since $\nu(p^{\gamma+1}) < +\infty$ and $\pi(\tilde{\delta}^{\gamma+1}) = \mathcal{M}_{\gamma+1} \leq \mathcal{M}_{\gamma_0+1} < +\infty$, it follows From Lebesgue's theorem that $\varepsilon$ is a continuous function on $\mathbb{R}$, which vanishes at $t = 0$. Point 4 is deduced from the above inequality and the equivalence of the norms $N_{\infty,\gamma}$ and $N_{1,\gamma}$. $\square$

Now assertions (a) and (c) of Proposition 7.1 follow directly from the results of Keller and Liverani (1999) which, moreover, assert that

$$\Pi_{(\gamma)}(t) = \frac{1}{2i\pi} \int_{\Gamma_1} (z - P(t))^{-1} \, dz$$

is a rank-1 bounded projection from $\mathcal{B}_\gamma$ onto $\mathrm{Ker}(P(t) - \lambda(t))$, and that $\nu(|\Pi_{(\gamma)}(t)1 - \Pi_{(\gamma)}(0)1|) = \nu(|\Pi_{(\gamma)}(t)1 - 1|)$ converges to 0 with $t$.

Therefore, for sufficiently small $|t|$, we have $\nu(\Pi_{(\gamma)}(t)1) \neq 0$, and we can set

$$v_{(\gamma)}(t) = \frac{1}{\nu(\Pi_{(\gamma)}(t)1)} \Pi_{(\gamma)}(t)1;$$

this function verifies condition (b) of Proposition 7.1. Assertion (d) and Corollary 7.2′ also follow from Keller and Liverani (1999). $\square$

We conclude this section with a result that will be useful for the study of the nonarithmeticity of $\xi$ (cf. Section 9).

PROPOSITION 7.4. *Assume that the conditions of Proposition 7.1 are satisfied and reinforced by $s + 1 < \gamma < \gamma_0$ and by the existence of $\tilde{\gamma}$, $\gamma < \tilde{\gamma} < \gamma_0$, such that $\mathcal{I}^{\tilde{\gamma}}(0, S) < +\infty$. Let $t \in \mathbb{R}$ be such that, for $P(t)$ acting on $\mathcal{B}_\gamma$, we have $r(P(t)) \geq 1$. Then $r(P(t)) = 1$ and $P(t)$ is quasi-compact on $\mathcal{B}_\gamma$.*

PROOF. Since $s + 1 < \gamma < \gamma_0$, we can suppose that $\tilde{\gamma}$ verifies $s + 1 + (\tilde{\gamma} - \gamma) \leq \gamma$. For convenience, we set $\tilde{N}(f) = N_{\infty,\gamma,\tilde{\gamma}}(f) = m_\gamma(f) + |f|_{\tilde{\gamma}}$ (Section 5.2).

The first inequality of Proposition 6.1, when applied to $\tilde{\gamma}$ and to the kernel $q(g, x) = e^{it\xi(g,x)}$, shows that

$$|P(t)f|_{\tilde{\gamma}} \leq \mathcal{I}^{\tilde{\gamma}}(0, 1)|f|_{\tilde{\gamma}},$$

with $\mathcal{I}^{\tilde{\gamma}}(0, 1) \leq \mathcal{I}^{\gamma_0}(0, 1) < +\infty$. Moreover, Proposition 6.7 applied to the couple $(\gamma, \tilde{\gamma}) = (\eta, \tilde{\eta})$ asserts that there exists a constant $C$ such that, for $t \in \mathbb{R}$ and $f \in \mathcal{B}_\gamma$, we have

$$m_\gamma(P(t)^{n_0}f) \leq \kappa_0^{n_0}m_\gamma(f) + C|t||f|_{\tilde{\gamma}}.$$

Setting $C' = C|t| + \mathcal{I}^{\tilde{\eta}}(0, 1)$, we get

$$\tilde{N}(P(t)^{n_0}f) \leq \kappa_0^{n_0}\tilde{N}(f) + C'|f|_{\tilde{\gamma}}.$$

From the fact that $P(t)$ is bounded on $(\mathcal{B}_\gamma, |\cdot|_{\tilde{\gamma}})$, and since the canonical embedding of $(\mathcal{B}_\gamma, \tilde{N})$ in $(\mathcal{B}_\gamma, |\cdot|_{\tilde{\gamma}})$ is compact (Lemma 5.4), we deduce by means of Corollary 1 of Hennion (1993) that, under the condition $r(P(t)) \geq 1$, $P(t)$ is quasi-compact on $\mathcal{B}_\gamma$, and that its essential spectral radius is $\leq \kappa_0$. Consequently, there exists an eigenvalue $\lambda$ of $P(t)$ such that $|\lambda| = r(P(t))$. Let $w \in \mathcal{B}_\gamma$ be an eigenfunction associated with $\lambda$. For $n \geq 1$, we have $|\lambda^n w| = |P(t)^n w| \leq P^n|w|$; hence $|\lambda^n| |w|_\gamma \leq |P^n|w| |_\gamma \leq \|P^n|w| \|_{\infty,\gamma}$. The spectral decomposition in Theorem 5.5 together with the equivalence of the considered norms on $\mathcal{B}_\gamma$ yield $\sup_n \|P^n|w| \|_{\infty,\gamma} < +\infty$. Hence $|\lambda| \leq 1$, and at last $r(P(t)) = 1$.  □

## 8. Taylor expansions for $v(\cdot)$, $\phi(\cdot)$, $N(\cdot)$.

The hypotheses in the subsequent statements will imply those of Proposition 7.1 and of its corollaries; thus, for small $|t|$, the eigenelements of the spectral decomposition described in Proposition 7.1 are defined. We are going to use the Taylor expansions of $P(\cdot)$ written in Proposition 6.3 to obtain the Taylor expansions for $v(\cdot)$, $\phi(\cdot)$ and $N(\cdot)$.

PROPOSITION 8.1 (First-order Taylor expansions). *Suppose that, for $\eta' < \eta$, the following condition holds:*

$$\mathcal{V}_1(\eta', \eta) : s + 1 \leq \eta' \leq \eta' + r < \eta \leq \gamma_0, \ \mathcal{I}^{\eta - r}(R, R + S) < +\infty.$$

*Then Proposition 7.1 applies to $P(t)$ acting on $\mathcal{B}_{\eta'}$, and the functions $v(\cdot)$, $\phi(\cdot)$ and $N(\cdot)$ from $I_{\eta'}$ in $(\mathcal{B}_{\eta'}, N_{\infty,\eta})$, $\mathcal{B}'_{\eta'}$ and $\mathcal{L}(\mathcal{B}_{\eta'}, \mathcal{B}_\eta)$, respectively, have a derivative at $t = 0$. Moreover, there exists a constant $K_1$ such that*

$$\forall n \geq 1, \ \forall t \in I_{\eta'} \qquad \|N(t)^n - N(0)^n\|_{\eta',\eta} \leq K_1|t|(\kappa'_0)^n.$$

PROPOSITION 8.2 (Second-order Taylor expansions). *Suppose that, for $\eta' < \eta$, the following condition holds:*

$$\mathcal{V}_2(\eta', \eta) : s + 1 \leq \eta' \leq \eta' + 2r < \eta \leq \gamma_0,$$

$$\mathcal{I}^{\eta - r}(R, R + S) + \mathcal{I}^{\eta'}(R^2, (R + S)R) < +\infty.$$

*Then Proposition 7.1 applies to $P(t)$ acting on $\mathcal{B}_{\eta'}$, and the functions $v(\cdot)$, $\phi(\cdot)$ and $N(\cdot)$, from $I_{\eta'}$ in $(\mathcal{B}_{\eta'}, N_{\infty, \eta})$, $\mathcal{B}'_{\eta'}$, and $\mathcal{L}(\mathcal{B}_{\eta'}, \mathcal{B}_\eta)$, respectively, have second-order Taylor expansions at $t = 0$. Moreover, we have, for all $t \in I_{\eta'}$ and $n \geq 1$,*

$$N(t)^n = N(0)^n + t\, N_{1,n} + \frac{t^2}{2} N_{2,n} + t^2 \varepsilon_n(t),$$

*with $N_{1,n}, N_{2,n}, \varepsilon_n(t) \in \mathcal{L}(\mathcal{B}_{\eta'}, \mathcal{B}_\eta)$, $\lim_{t \to 0} \sup_{n \geq 1} \|\varepsilon_n(t)\|_{\eta', \eta} = 0$, and $\sup_{n \geq 1} \|N_{j,n}\|_{\eta', \eta} < +\infty$ for $j = 1, 2$.*

The rest of this section is devoted to the proofs of these propositions. Recall that $\mathcal{R} = \{z : z \in \mathbb{C}, |z| \geq \kappa'_0, |z - 1| \geq 1 - \kappa''_0\}$. For $\gamma \in ]0, \gamma_0]$, we set

$$J_\gamma = \sup_{z \in \mathcal{R}} \|(z - P)^{-1}\|_\gamma < +\infty \qquad \text{(Theorem 5.5)}.$$

Under condition $\mathcal{V}_1(\eta', \eta)$ or $\mathcal{V}_2(\eta', \eta)$, we have $\mathcal{I}^{\eta'}(0, S) \leq \mathcal{I}^{\eta - r}(0, S) < +\infty$ and $s + 1 \leq \eta' \leq \gamma_0$. Consequently, Proposition 7.1 applies to $P(t)$ acting on $\mathcal{B}_{\eta'}$. In particular, for $t \in I_{\eta'}$ and for $z \in \mathcal{R}$, $(z - P(t))$ is invertible on $\mathcal{B}_{\eta'}$, and we have

$$M_{\eta'} = \sup\{\|(z - P(t))^{-1}\|_{\eta'}, \ t \in I_{\eta'}, \ z \in \mathcal{R}\} < +\infty.$$

We shall need the following formula. Let $\mathcal{B}$ be a Banach space. If $U$ and $V$ are bounded operators on $\mathcal{B}$ such that $U$ and $U - V$ are invertible, we have

$$(*) \qquad (U - V)^{-1} = \sum_{k=0}^{n} (U^{-1}V)^k U^{-1} + (U^{-1}V)^{n+1}(U - V)^{-1}.$$

Actually, if $W \in \mathcal{L}(\mathcal{B})$, we have $I - W^{n+1} = \sum_{k=0}^{n} W^k (I - W)$, and hence, if $I - W$ is invertible,

$$(I - W)^{-1} = \sum_{k=0}^{n} W^k + W^{n+1}(I - W)^{-1}.$$

The claimed formula follows from the relation $(U - V)^{-1} = (I - U^{-1}V)^{-1}U^{-1}$ and the above equality.

In the proofs below, we shall apply $(*)$ with $U = z - P$, $V = P(t) - P$, and thus $U - V = z - P(t)$. Observe that, in the sequel, all the space parameters $\gamma$ are between $s+1$ and $\gamma_0$, so that conditions $\mathcal{U}_0(\eta',\eta)$ and $\mathcal{U}_n(\eta',\eta)$, $n \geq 1$, of Propositions 6.2 and 6.3 can be rewritten as

$$\mathcal{U}_0(\eta',\eta): \eta' < \eta, \ \mathcal{I}^{\eta'}(0,S) < +\infty,$$

$$\mathcal{U}_n(\eta',\eta): \eta' + nr < \eta, \ \mathcal{I}^{\eta'}(R^n,(R+S)R^{n-1}) < +\infty.$$

Otherwise notice that, if $\eta' < \eta_1 < \eta$ and if $T \in \mathcal{L}(\mathcal{B}_{\eta_1},\mathcal{B}_\eta)$, then $T \in \mathcal{L}(\mathcal{B}_{\eta'},\mathcal{B}_\eta)$ and $\|T\|_{\eta',\eta} \leq \|T\|_{\eta_1,\eta}$.

PROOF OF PROPOSITION 8.1.    The next lemma gives a first-order Taylor expansion for the resolvent $(z-P(t))^{-1}$. We set $R(z,t) = (z-P(t))^{-1}$ and $R(z) = R(z,0) = (z-P)^{-1}$.

LEMMA 8.3.    *Under condition $\mathcal{V}_1(\eta',\eta)$, there exists a continuous function $R'_z$ from $\mathcal{R}$ to $\mathcal{L}(\mathcal{B}_{\eta'},\mathcal{B}_\eta)$, such that we have*

$$\lim_{t \to 0} \frac{1}{|t|} \sup_{z \in \mathcal{R}} \|(z-P(t))^{-1} - (z-P)^{-1} - tR'_z\|_{\eta',\eta} = 0.$$

PROOF.    Setting $n$ to 1 and $U$ and $V$ to the values indicated a few lines above, the formula $(*)$ gives, for $z \in \mathcal{R}$ and $t \in I_{\eta'}$,

$$R(z,t) = R(z) + R(z)(P(t)-P)R(z) + R(z)(P(t)-P)R(z)(P(t)-P)R(z,t).$$

As, by assumption, $\eta' + r < \eta$, we can choose $\eta_1$ such that $\eta' < \eta_1 \leq \eta_1 + r < \eta$.

Condition $\mathcal{U}_0(\eta',\eta_1)$ is verified because $\eta' < \eta_1$ and $\mathcal{I}^{\eta'}(0,S) \leq \mathcal{I}^{\eta-r}(0,S) < +\infty$; hence $\lim_{t\to 0} \|P(t) - P\|_{\eta',\eta_1} = 0$. Condition $\mathcal{U}_1(\eta_1,\eta)$ holds because $\eta_1 + r < \eta$ and $\mathcal{I}^{\eta_1}(R,R+S) \leq \mathcal{I}^{\eta-r}(R,R+S) < +\infty$; hence $P(t) - P = tL_1 + \Upsilon_1(t)$, with $L_1, \Upsilon_1(t) \in \mathcal{L}(\mathcal{B}_{\eta_1},\mathcal{B}_\eta)$ and $\lim_{t\to 0} |t|^{-1}\|\Upsilon_1(t)\|_{\eta_1,\eta} = 0$ (Proposition 6.3). Now we write

$$R(z,t) = R(z) + tR'_z + \Theta_1(z,t) + \Theta_2(z,t),$$

with $R'_z = R(z)L_1R(z)$, and

$$\Theta_1(z,t) = R(z)\Upsilon_1(t)R(z),$$

$$\Theta_2(z,t) = R(z)(P(t)-P)R(z)(P(t)-P)R(z,t).$$

Since $L_1 \in \mathcal{L}(\mathcal{B}_{\eta_1},\mathcal{B}_\eta) \subset \mathcal{L}(\mathcal{B}_{\eta'},\mathcal{B}_\eta)$ and since $R(\cdot)$ is continuous from $\mathcal{R}$ to both $\mathcal{L}(\mathcal{B}_{\eta'})$ and $\mathcal{L}(\mathcal{B}_\eta)$, $R'_z$ is continuous from $\mathcal{R}$ to $\mathcal{L}(\mathcal{B}_{\eta'},\mathcal{B}_\eta)$. For $t \in I_{\eta'}$, $0 < |t| \leq 1$, and $z \in \mathcal{R}$, we have

$$|t|^{-1}\|\Theta_1(z,t)\|_{\eta',\eta} \leq |t|^{-1}\|\Theta_1(z,t)\|_{\eta_1,\eta} \leq J_\eta \ |t|^{-1}\|\Upsilon_1(t)\|_{\eta_1,\eta}J_{\eta_1},$$

$$|t|^{-1}\|\Theta_2(z,t)\|_{\eta',\eta} \leq J_\eta(\|L_1\|_{\eta_1,\eta} + |t|^{-1}\|\Upsilon_1(t)\|_{\eta_1,\eta})J_{\eta_1}\|P(t) - P\|_{\eta',\eta_1}M_{\eta'}.$$

The second members do not depend on $z \in \mathcal{R}$ and converge to 0 with $t$, this proves the lemma.    $\square$

To establish Proposition 8.1, we now use the formulae of Corollary $7.2'$. More precisely, the linear maps $\Pi(t)$ and $N(t)$ of the corollary are considered here as elements of $\mathcal{L}(\mathcal{B}_{\eta'}, \mathcal{B}_{\eta})$ since they may be viewed as integrals of functions with values in $\mathcal{L}(\mathcal{B}_{\eta'}, \mathcal{B}_{\eta})$.

Then Lemma 8.3 shows that $\Pi(\cdot)$ has a derivative at $t = 0$ as an $\mathcal{L}(\mathcal{B}_{\eta'}, \mathcal{B}_{\eta})$-valued function. Thus $\Pi(\cdot)^*$ has a derivative at $t = 0$ as an $\mathcal{L}(\mathcal{B}'_\eta, \mathcal{B}'_{\eta'})$-valued function. This proves the first-order Taylor expansions of $v(\cdot)$ and $\phi(\cdot)$. The existence of a derivative for $N(\cdot)$ at $t = 0$ follows in a similar way from Lemma 8.3. On the other hand, from the integral formula $N(t)^n = \frac{1}{2i\pi} \int_{\Gamma_0} z^n (z - P(t))^{-1} \, dz$, we deduce the existence of a constant $K$ such that, for $n \geq 1$ and $t \in I_{\eta'}$,

$$(\kappa_0')^{-(n+1)} \|N(t)^n - N(0)^n\|_{\eta',\eta} \leq \sup_{z \in \Gamma_0} \|(z - P(t))^{-1} - (z - P)^{-1}\|_{\eta',\eta}$$

$$\leq |t| \left( \sup_{z \in \Gamma_0} \|R'_z\|_{\eta',\eta} + K \right);$$

hence the inequality of Proposition 8.1.   $\square$

PROOF OF PROPOSITION 8.2.   As above, we start with a Taylor expansion of the resolvent $(z - P(t))^{-1}$.

LEMMA 8.4.   *Under condition* $\mathcal{V}_2(\eta', \eta)$, *there exist continuous functions* $R'_z$ *and* $R''_z$ *from* $\mathcal{R}$ *to* $\mathcal{L}(\mathcal{B}_{\eta'}, \mathcal{B}_{\eta})$, *such that we have*

$$\lim_{t \to 0} \frac{1}{t^2} \sup_{z \in \mathcal{R}} \left\| (z - P(t))^{-1} - (z - P)^{-1} - tR'_z - \frac{t^2}{2} R''_z \right\|_{\eta',\eta} = 0.$$

PROOF.   Retaining the notation of Lemma 8.3 but setting $n$ to 2, the formula $(*)$ gives, for $z \in \mathcal{R}$ and $t \in I_{\eta'}$,

$$R(z,t) = R(z) + R(z)(P(t) - P)R(z) + R(z)(P(t) - P)R(z)(P(t) - P)R(z)$$
$$+ R(z)(P(t) - P)R(z)(P(t) - P)R(z)(P(t) - P)R(z,t).$$

Since $\eta' + 2r < \eta$, we can choose $\eta_1$ and $\eta_2$ such that $\eta' < \eta_1 \leq \eta_1 + r < \eta_2 \leq \eta_2 + r < \eta$.

The condition $\mathcal{U}_2(\eta', \eta)$ is verified; hence by Proposition 6.3,

$$P(t) - P = tL_1 + \frac{t^2}{2} L_2 + \Upsilon_2(t),$$

with $P, P(t), L_1, L_2, \Upsilon_2(t) \in \mathcal{L}(\mathcal{B}_{\eta'}, \mathcal{B}_{\eta})$ and $\lim_{t \to 0} t^{-2} \|\Upsilon_2(t)\|_{\eta',\eta} = 0$.

The conditions $\mathcal{U}_1(\eta_1, \eta_2)$ and $\mathcal{U}_1(\eta_2, \eta)$ are satisfied since we have $\eta_1 + r < \eta_2$, $\eta_2 + r < \eta$ and $\mathcal{I}^{\eta_1}(R, R + S) \leq \mathcal{I}^{\eta_2}(R, R + S) \leq \mathcal{I}^{\eta - r}(R, R + S)$. Then Proposition 6.3 with $n = 1$ shows that

$$P(t) - P = tL_1 + \Upsilon_1(t),$$

with $L_1, \Upsilon_1(t) \in \mathcal{L}(\mathcal{B}_{\eta_1}, \mathcal{B}_{\eta_2}) \cap \mathcal{L}(\mathcal{B}_{\eta_2}, \mathcal{B}_{\eta})$ and $\lim_{t \to 0} |t|^{-1} \|\Upsilon_1(t)\|_{\eta_1, \eta_2} = \lim_{t \to 0} |t|^{-1} \|\Upsilon_1(t)\|_{\eta_2, \eta} = 0$.

At last, since $\mathcal{I}^{\eta'}(0, S) \le \mathcal{I}^{\eta-r}(0, S) < +\infty$, the condition $\mathcal{U}_0(\eta', \eta_1)$ holds and Proposition 6.2 ensures that

$$\lim_{t \to 0} \|P(t) - P\|_{\eta', \eta_1} = 0.$$

We get

$$\begin{aligned}
R(z, t) = R(z) &+ R(z)\Big(tL_1 + \frac{t^2}{2}L_2 + \Upsilon_2(t)\Big)R(z) \\
&+ R(z)(tL_1 + \Upsilon_1(t))R(z)(tL_1 + \Upsilon_1(t))R(z) \\
&+ R(z)(tL_1 + \Upsilon_1(t))R(z)(tL_1 + \Upsilon_1(t))R(z)(P(t) - P)R(z, t),
\end{aligned}$$

hence

$$R(z, t) = R(z) + tR'_z + \frac{t^2}{2}R''_z + \sum_{k=1}^{5}\Theta_k(z, t),$$

with

$$R'_z = R(z)L_1 R(z), \qquad R''_z = R(z)L_2 R(z) + 2R(z)L_1 R(z)L_1 R(z),$$

and

$$\begin{aligned}
\Theta_1(z, t) &= R(z)\Upsilon_2(t)R(z), \\
\Theta_2(z, t) &= tR(z)L_1 R(z)\Upsilon_1(t)R(z), \\
\Theta_3(z, t) &= tR(z)\Upsilon_1(t)R(z)L_1 R(z), \\
\Theta_4(z, t) &= R(z)\Upsilon_1(t)R(z)\Upsilon_1(t)R(z), \\
\Theta_5(z, t) &= R(z)(tL_1 + \Upsilon_1(t))R(z)(tL_1 + \Upsilon_1(t))R(z)(P(t) - P)R(z, t).
\end{aligned}$$

Since $L_1 \in \mathcal{L}(\mathcal{B}_{\eta'}, \mathcal{B}_{\eta}) \cap \mathcal{L}(\mathcal{B}_{\eta'}, \mathcal{B}_{\eta_2}) \cap \mathcal{L}(\mathcal{B}_{\eta_2}, \mathcal{B}_{\eta})$, $L_2 \in \mathcal{L}(\mathcal{B}_{\eta'}, \mathcal{B}_{\eta})$ and $R(\cdot)$ is continuous from $\mathcal{R}$ to $\mathcal{L}(\mathcal{B}_{\eta'})$, $\mathcal{L}(\mathcal{B}_{\eta_2})$, and $\mathcal{L}(\mathcal{B}_{\eta})$, the functions $R'_z$ and $R''_z$ are continuous from $\mathcal{R}$ to $\mathcal{L}(\mathcal{B}_{\eta'}, \mathcal{B}_{\eta})$.

We have, for $t \in I$, $0 < |t| \le 1$, and $z \in \mathcal{R}$

$$\begin{aligned}
t^{-2}\|\Theta_1(z, t)\|_{\eta', \eta} &\le J_\eta(t^{-2}\|\Upsilon_2(t)\|_{\eta', \eta})J_{\eta'}, \\
t^{-2}\|\Theta_2(z, t)\|_{\eta', \eta} &\le t^{-2}\|\Theta_2(z, t)\|_{\eta_1, \eta} \\
&\le J_\eta\|L_1\|_{\eta_2, \eta}J_{\eta_2}(|t|^{-1}\|\Upsilon_1(t)\|_{\eta_1, \eta_2})J_{\eta_1}, \\
t^{-2}\|\Theta_3(z, t)\|_{\eta', \eta} &\le t^{-2}\|\Theta_3(z, t)\|_{\eta_1, \eta} \\
&\le J_\eta(|t|^{-1}\|\Upsilon_1(t)\|_{\eta_2, \eta})J_{\eta_2}\|L_1\|_{\eta_1, \eta_2}J_{\eta_1}, \\
t^{-2}\|\Theta_4(z, t)\|_{\eta', \eta} &\le t^{-2}\|\Theta_4(z, t)\|_{\eta_1, \eta} \\
&\le J_\eta(|t|^{-1}\|\Upsilon_1(t)\|_{\eta_2, \eta})J_{\eta_2}(|t|^{-1}\|\Upsilon_1(t)\|_{\eta_1, \eta_2})J_{\eta_1}, \\
t^{-2}\|\Theta_5(z, t)\|_{\eta', \eta} &\le J_\eta K_{\eta_2, \eta}J_{\eta_2}K_{\eta_1, \eta_2}J_{\eta_1}\|P(t) - P\|_{\eta', \eta_1}M_{\eta'},
\end{aligned}$$

with $K_{a,b} = \sup\{\|L_1\|_{a,b} + |t|^{-1}\|\Upsilon_1(t)\|_{a,b}, t \in I, |t| \leq 1\}$.

This proves the lemma because the right-hand members do not depend on $z \in \mathcal{R}$ and tend to 0 with $t$. $\quad\square$

Let us now complete the proof of Proposition 8.2. The Taylor expansions for $v(\cdot)$, $\phi(\cdot)$ and $N(\cdot)$ can be deduced from the formulae of Corollary 7.2′. We just specify how to get the expansion for $N(\cdot)^n$. Using integration in $\mathcal{L}(\mathcal{B}_{\eta'}, \mathcal{B}_{\eta})$, we set

$$N_{1,n} = \frac{1}{2i\pi}\int_{\Gamma_0} z^n R_z'\, dz \quad \text{and} \quad N_{2,n} = \frac{1}{2i\pi}\int_{\Gamma_0} z^n R_z''\, dz.$$

We have $\|N_{1,n}\|_{\eta',\eta} \leq \kappa_0'^{n+1}\sup_{z\in\Gamma_0}\|R_z'\|_{\eta',\eta}$ and $\|N_{2,n}\|_{\eta',\eta} \leq \kappa_0'^{n+1}\times\sup_{z\in\Gamma_0}\|R_z''\|_{\eta',\eta}$. Lemma 8.4 yields

$$\|\varepsilon_n(t)\|_{\eta',\eta} = \frac{1}{2\pi t^2}\left\|\int_{\Gamma_0} z^n\left((z-P(t))^{-1} - (z-P)^{-1} - tR_z' - \frac{t^2}{2}R_z''\right)dz\right\|_{\eta',\eta}$$

$$\leq \frac{\kappa_0'^{n+1}}{t^2}\sup_{z\in\mathcal{R}}\left\|(z-P(t))^{-1} - (z-P)^{-1} - t\,R_z' - \frac{t^2}{2}R_z''\right\|_{\eta',\eta}.$$

Since $\kappa_0'^{n+1} \leq 1$, we conclude that $\lim_{t\to 0}\sup_{n\geq 1}\|\varepsilon_n(t)\|_{\eta',\eta} = 0$. $\quad\square$

**9. Extensions and proofs of Theorems A, B, C, S.** We return to the context of Sections 1 and 2. Theorems A′, B′, C′ below concern the behaviour of the sequence of r.v.'s $((Z_n, S_n^Z))_n$.

9.1. *Theorems* A′, B′, C′, S′. Neglecting the technical parameter $\lambda_0$ of the preceding sections, we may define $\mathcal{B}_\gamma$ as the space of locally Lipschitz $\mathbb{C}$-valued functions $f$ on $M$ such that

$$\ell_\gamma(f) = \sup\left\{\frac{|f(x) - f(y)|}{d(x,y)(1 + d(x,x_0))^\gamma(1 + d(y,x_0))^\gamma}, x,y \in M, x \neq y\right\} < +\infty,$$

endowed with the norm

$$\|f\|_{\infty,\gamma} = \ell_\gamma(f) + \sup_{x\in M}\frac{|f(x)|}{(1 + d(x,x_0))^{\gamma+1}};$$

this norm is clearly equivalent to the ones previously defined on $\mathcal{B}_\gamma$.

Recall that, we set $\tilde{\delta}(g) = 1 + c(g) + d(gx_0, x_0)$. As previously we can omit $\lambda_0$ in the definition of the numbers $\mathcal{I}^\tau(U, V)$ (Section 6.2) by replacing now the function $\delta$ by $\tilde{\delta}$, that is, by replacing $\mathcal{I}^\tau(U, V)$ by $\mathcal{J}^\tau(U, V)$, already used in Section 2, and defined by

$$\mathcal{J}^\tau(U, V) = \int_G U(g)c(g)\tilde{\delta}(g)^{2\tau}\, d\pi(g) + \int_G V(g)\tilde{\delta}(g)^{\tau+1}\, d\pi(g).$$

If $(V, \|\cdot\|)$ is a normed linear space and if $\alpha > 0$, we shall denote by $V(\alpha)$ the closed ball in $V$ with radius $\alpha$ centered at 0. We name $\mathcal{C}_{\downarrow 2}(\mathbb{R})$ the space of $\mathbb{C}$-valued continuous functions $h$ on $\mathbb{R}$ such that $\lim_{|u| \to +\infty} u^2 h(u) = 0$.

Under the hypotheses of the next statements, the real number $m = \int_M \int_G \xi(g,x) \, d\pi(g) \, d\nu(x)$ is defined, and supposed to be zero.

Recall that Condition $\mathcal{H}(\gamma_0)$ holds if there exist $\gamma_0 \in \mathbb{R}_+^*$ and $n_0 \in \mathbb{N}^*$ such that

$$\mathcal{M}_{\gamma_0 + 1} = \pi(\tilde{\delta}^{\gamma_0 + 1}) < +\infty,$$

$$\mathcal{M}'_{2\gamma_0 + 1} = \pi(c\tilde{\delta}^{2\gamma_0}) < +\infty,$$

$$\mathcal{C}^{(n_0)}_{2\gamma_0 + 1} = \pi^{*n_0}(c \max\{c, 1\}^{2\gamma_0}) < 1.$$

THEOREM A$'$ (Central limit).  *Assume $\mathcal{H}(\gamma_0)$ with $\gamma_0 > r + \max\{r, s+1\}$ and that*

$$\int_G R^2 \, d\pi < +\infty, \qquad \mathcal{J}^{\gamma_0 - r}(R, R+S) < +\infty.$$

*Then, there exists $\sigma^2 \geq 0$ such that, if the r.v. $Z$ satisfies $\mathbb{E}[d(Z, x_0)^{\gamma_0 + 1}] < +\infty$, we have, for $f \geq 0$, $f \in \bigcup_{\gamma < \gamma_0 - r} \mathcal{B}_\gamma$, and for any bounded continuous function $h$ on $\mathbb{R}$,*

$$\lim_n \mathbb{E}\Big[ f(Z_n) h\Big(\frac{S_n^Z}{\sqrt{n}}\Big) \Big] = \nu(f)\mathcal{N}(0, \sigma^2)(h).$$

*If $h \in C_{\downarrow 2}(\mathbb{R})$, this convergence holds uniformly when $(\mu, f)$ ranges over $\mathcal{B}'_{\gamma_0}(\alpha) \times \mathcal{B}_\gamma(\alpha)$.*

THEOREM B$'$ (Central limit with a rate of convergence).  *Assume $\mathcal{H}(\gamma_0)$ with $\gamma_0 > 3r + \max\{r, s+1\}$ and that*

$$\int_G R^3 \, d\pi < +\infty, \qquad \mathcal{J}^{\gamma_0 - r}(R, R+S) + \mathcal{J}^{\gamma_0 - 2r}(R^2, (R+S)R) < +\infty.$$

*Then, if $\sigma^2 > 0$, the assertion of Theorem B holds.*

*Moreover, if $Z$ has the distribution $\nu$, then, for $0 < \gamma < \gamma_0 - r$, there exists a positive constant $C_\gamma$ such that, for $f \in \mathcal{B}_\gamma$, $f \geq 0$, satisfying $\nu(f) > 0$, we have*

$$\sup_{u \in \mathbb{R}} |\mathbb{E}[f(Z_n)\mathbb{1}_{[S_n^Z \leq u\sigma\sqrt{n}]}] - \nu(f)\mathcal{N}(0,1)(]-\infty, u])| \leq \frac{C_\gamma \|f\|_{\infty, \gamma}}{\sqrt{n}}.$$

The statement of the local limit theorem appeals to the nonarithmeticity condition for $\xi$ with respect to the space $\mathcal{B}_\gamma$ for $\gamma \in ]s+1, \gamma_0 - r[$:

CONDITION $(N-A)_\gamma$.  *There is no $t \in \mathbb{R} \setminus \{0\}$, no $\lambda \in \mathbb{C}$, $|\lambda| = 1$, and no bounded function $w$ in $\mathcal{B}_\gamma$ with nonzero constant modulus on the support $\Sigma_\nu$ of $\nu$, such that we have, for all $x \in \Sigma_\nu$ and for all $n \geq 1$,*

$$e^{itS_n^x}w(R_n x) = \lambda^n w(x), \qquad \mathbb{P}\text{-}a.s.$$

THEOREM C$'$ (Local central limit).  *Assume that the hypotheses of Theorem A$'$ are satisfied. Let $\gamma$ be a real number verifying $\max\{r, s+1\} < \gamma < \gamma_0 - r$ and such that Condition $(N$–$A)_\gamma$ is fulfilled.*

*If $\sigma^2 > 0$, and if $Z$ is such that $\mathbb{E}[d(Z, x_0)^{\gamma_0+1}] < +\infty$, then for all $f \geq 0$, $f \in \mathcal{B}_\gamma$, and for all $h \in \mathcal{C}_{\downarrow 2}(\mathbb{R})$, we have*

$$\limsup_n \sup_{u \in \mathbb{R}} |\sigma\sqrt{2\pi n}\mathbb{E}[f(Z_n)h(S_n^Z - u)] - e^{-u^2/(2n\sigma^2)}\nu(f)\mathcal{L}(h)| = 0,$$

*and this convergence holds uniformly when $(\mu, f)$ ranges over $\mathcal{B}'_{\gamma_0}(\alpha) \times \mathcal{B}_\gamma(\alpha)$.*

THEOREM S$'$.  *Assume $\mathcal{H}(\gamma_0)$ with $\gamma_0 > 2r + s + 1$ and that*

$$\mathcal{J}^{\gamma_0-r}(R, R+S) + \mathcal{J}^{\gamma_0-2r}(R^2, (R+S)R) < +\infty.$$

*Then the assertions of Theorem S hold with $\tilde{\xi}_1 \in \mathcal{B}_{\gamma_0-r}$ in point* (i).

9.2.  *Proofs of Theorems* A$'$, B$'$, C$'$.  These proofs are based on expansions of the characteristic function of the r.v $S_n^Z$.

PROPOSITION 9.1.  1.  *Assume that the hypotheses of Theorem A$'$ are fulfilled. Let the parameter $\gamma$ verify $\max\{r, s+1\} < \gamma < \gamma_0 - r$.*

*Then there exist an open interval $I_\gamma$ containing $t = 0$, a $\mathbb{C}$-valued function $\lambda(\cdot)$ and $\mathcal{L}(\mathcal{B}_\gamma)$-valued functions $L(\cdot)$, $N(\cdot)$, defined on this interval, such that, if the distribution $\mu$ of the r.v. $Z$ verifies $\mu(d(\cdot, x_0)^{\gamma_0+1}) < +\infty$, we have, for $n \geq 1$, $t \in I_\gamma$, and for $f \in \mathcal{B}_\gamma$,*

$$\mathbb{E}[f(Z_n)e^{itS_n^Z}] = \langle \mu, P(t)^n f \rangle = \lambda(t)^n(\nu(f) + \langle \mu, L(t)f \rangle) + \langle \mu, N(t)^n f \rangle.$$

*For all $t \in I_\gamma$, we have $|\lambda(t)| \leq 1$; there exists a real positive number $\sigma^2 \geq m^2$ such that*

$$\lambda(t) = 1 + imt - \sigma^2\frac{t^2}{2} + o(t^2),$$

*and there exists a positive constant $c_\gamma$ such that:*

(i) *if, either $f = 1$ and $\mu \in \mathcal{B}'_{\gamma_0}$, or $f \in \mathcal{B}_\gamma$ and $\mu = \nu$, then*

$$|\langle \mu, N(t)^n f \rangle| \leq c_\gamma(\kappa'_0)^n \inf\{|t|, 1\}\|\mu\|_{\infty,\gamma_0}\|f\|_{\infty,\gamma},$$

(ii) $\|N(t)^n\|_\gamma \leq c_\gamma(\kappa'_0)^n$,

(iii) $\|L(t)\|_{\gamma,\gamma_0} \leq c_\gamma \inf\{|t|, 1\}.$

*Moreover, if $m = 0$ and $\sigma^2 > 0$, then, for any real number $t$ such that $\frac{t}{\sigma} \in I_\gamma$, we have:*

(iv) $|\lambda(\frac{t}{\sigma})| \leq e^{-t^2/4}.$

2. *Suppose that the hypotheses of Theorem* B′ *hold. Then, if $m = 0$ and $\sigma^2 > 0$, there exists a constant $C_1$ such that we have, for all real $t$ such that $\frac{t}{\sigma\sqrt{n}} \in I_\gamma$,*

(v) $|\lambda(\frac{t}{\sigma\sqrt{n}})^n - e^{-t^2/2}| \leq \frac{C_1}{\sqrt{n}}|t|^3 e^{-t^2/4}.$

Assume this proposition for a while. To prove Theorems A′, B′, C′, we have only to use the method of Hennion and Hervé [[2001](), Section IV.2 and Chapter VI], which is an adaptation of standard Fourier techniques for sums of i.i.d. r.v.'s. As already mentioned in Section 4.2, we consider here the Fourier kernels $P(t)$ instead of the Fourier kernels $Q(t)$ associated with $\xi$ and the probability transition $Q$ on $G \times M$ defined in Section 4.2. Yet the needed changes are obvious, and we shall not develop the argumentation; we only specify some points.

First, the distribution $\mu$ of $Z$ defines an element of $\mathcal{B}'_{\gamma_0}$ if and only if $\mathbb{E}[d(Z, x_0)^{\gamma_0+1}] < +\infty$, and, in this case, $\|\mu\|_{\infty,\gamma_0} = \mathbb{E}[(1 + d(Z, x_0))^{\gamma_0+1}]$. Actually, we have, for $f \in \mathcal{B}_{\gamma_0}$, $|f| \leq \|f\|_{\infty,\gamma_0}(1 + d(\cdot, x_0))^{\gamma_0+1}$; hence $\mu(|f|) \leq \|f\|_{\infty,\gamma_0}\mathbb{E}[(1 + d(Z, x_0))^{\gamma_0+1}]$.

Second, because of the topological embedding of the spaces $\mathcal{B}_\gamma$, in the proofs of Theorems A′ and B′, it will be sufficient to consider the case where the function $f$ is in a space $\mathcal{B}_\gamma$ with $\gamma \in ]\max\{r, s+1\}, \gamma_0 - r[$.

At last, in the proof of Theorem C′, it is necessary to have some control on the behavior of $P(t)$, for all $t \in \mathbb{R}$. The following lemma shows how this is related to Condition (N–A)$_\gamma$.

LEMMA 9.1′. *Assume conditions of Theorem* C′ *except Condition* (N–A)$_\gamma$. *Then $P(t)$ is a bounded operator of $\mathcal{B}_\gamma$ for all $t \in \mathbb{R}$. Let $t \in \mathbb{R}$ such that $r(P(t)) \geq 1$. Then there exist $\lambda \in \mathbb{C}$, $|\lambda| = 1$ and a bounded function $w \in \mathcal{B}_\gamma$, with nonzero constant modulus on the support $\Sigma_\nu$ of $\nu$, such that we have, for all $x \in \Sigma_\nu$ and all $n \geq 1$,*

$$e^{itS_n^x}w(R_n x) = \lambda^n w(x), \qquad \mathbb{P}\text{-}a.s.$$

*Consequently, under Condition* (N–A)$_\gamma$, *for all $t \in \mathbb{R} \setminus \{0\}$, we have $r(P(t)) < 1$.*

PROOF. Let $\gamma \in ]s+1, \gamma_0 - r[$. The inequality $\mathcal{J}^\gamma(0, S) \leq \mathcal{J}^{\gamma_0-r}(0, S) < +\infty$ together with Proposition 7.1 shows that the Fourier kernels $P(t)$ act continuously on $\mathcal{B}_\gamma$ for all $t \in \mathbb{R}$.

By Proposition 7.4, if $r(P(t)) \geq 1$, then $r(P(t)) = 1$, and $P(t)$ is quasi-compact. Consequently, there exist $w \in \mathcal{B}_\gamma \setminus \{0\}$, and $\lambda \in \mathbb{C}$, $|\lambda| = 1$, such that, for all $n \geq 1$, we have

$$P(t)^n w = \lambda^n w.$$

It follows that $|w| \leq P^n |w|$. Since, by Theorem 5.1, the sequence $(P^n |w|)_{n \geq 1}$ converges pointwise to $\nu(|w|)$, we get $|w| \leq \nu(|w|)$, so that $w$ is bounded. From the above and equality $\nu(\nu(|w|)1_M - |w|) = 0$, we deduce that $\nu(\{x : x \in M, |w(x)| = \nu(|w|)\}) = 1$; thus $|w|$ is a nonzero constant function on $\Sigma_\nu$. For $x \in \Sigma_\nu$ and $n \geq 1$, we write

$$\mathbb{E}\left[1 - \frac{e^{itS_n^x}w(R_n x)}{\lambda^n w(x)}\right] = 1 - \frac{P(t)^n w(x)}{\lambda^n w(x)} = 0.$$

Since $\left|\frac{e^{itS_n^x}w(R_n x)}{\lambda^n w(x)}\right| = 1$, it follows that $e^{itS_n^x}w(R_n x) = \lambda^n w(x)$, $\mathbb{P}$-a.s. $\quad\square$

To be complete on the properties required for local theorem, one needs to establish the following.

LEMMA 9.1″. *Under the conditions of Theorem* C′, *for every compact subset $K$ of $\mathbb{R}^*$:*

(i) *We have $r_K = \sup\{r(P(t)), t \in K\} < 1$.*

(ii) *There exists $C \geq 0$ and $\rho_K < 1$ such that we have, for all $n \geq 1$,* $\sup_{t \in K} \|P(t)^n\| \leq C\rho_K^n$.

PROOF. (i) Suppose that $\sup_{t \in K} r(P(t)) \geq 1$. Then, by Lemma 9.1′, $\sup_{t \in K} r(P(t)) = 1$, thus there exists a sequence $(\tau_k)_k$ in $K$ such that $\lim_k r(P(\tau_k)) = 1$. For each $k \geq 1$ consider a spectral value $\lambda_k$ of $P(\tau_k)$ satisfying $|\lambda_k| = r(P(\tau_k))$. By compactness, one can suppose that $(\tau_k)_k$ and $(\lambda_k)_k$ converge. Set $t_0 = \lim_k \tau_k$, $\lambda = \lim_k \lambda_k$, and observe that $t_0 \in K$, thus $t_0 \neq 0$, and $|\lambda| = 1$.

We are going to show that the perturbation theorem of Keller and Liverani (1999) applies to the action, on a certain space $\mathcal{B}_\gamma$, of the family $\{P(t), t \in \mathbb{R}\}$ when $t \to t_0$. It will follows from this result, see page 145 of the above cited paper, that $\lambda$ is a spectral value of $P(t_0)$. But since $t_0 \neq 0$ and $|\lambda| = 1$, this will contradict Lemma 9.1′, so we shall get point (i).

Let $\gamma, \tilde{\gamma}$ be such that $s + 1 + (\tilde{\gamma} - \gamma) \leq \gamma < \tilde{\gamma} < \gamma_0 - r$. We establish that $\{P(t), t \in \mathbb{R}\}$ acting on $\mathcal{B}_\gamma$ satisfies the four assertions of Lemma 7.3, where $0$ is replaced by $t_0 \in \mathbb{R}^*$, and the norm $N_{1,\gamma}(\cdot)$ [resp. $\nu(\cdot)$] is replaced by $N_{\infty,\gamma,\tilde{\gamma}}$ (resp. $|\cdot|_{\tilde{\gamma}}$).

1. Using the inequality $|P(t)^n f| \leq P^n |f| \leq |f|_{\tilde{\gamma}} P^n(p^{\tilde{\gamma}+1})$ and assertion (c) of Lemma 4.2 (observe that $p^{\tilde{\gamma}}$ and $1 + \phi_\lambda$ are equivalent), one easily proves that $\sup_{n \geq 1} |P^n(p^{\tilde{\gamma}+1})|_{\tilde{\gamma}} < +\infty$. It follows that $\{P(t)^n, t \in \mathbb{R}, n \geq 1\}$ is uniformly bounded on $(\mathcal{B}_\gamma, |\cdot|_{\tilde{\gamma}})$.

2. Proposition 6.7 implies the second point of Lemma 7.3 [with $|\cdot|_{\tilde{\gamma}}$ instead of $\nu(|\cdot|)$].

3. If $r(P(t)) > \vartheta_0^{1/n_0}$, where $\vartheta_0 < 1$ is the real number in Proposition 6.7, it follows from Lemma 5.4, from the preceding assertion, and from Hennion (1993), that the essential spectral radius of $P(t)$ is $\leq \vartheta_0^{1/n_0}$. If $r(P(t)) \leq \vartheta_0^{1/n_0}$, this is also valid because the essential spectral radius is always less than the spectral radius.

4. In the same way as Proposition 6.2, it can be proved that there exists a real continuous function $\varepsilon(\cdot)$, vanishing at $t = t_0$, such that we have $\|P(t)f - P(t_0)f\|_{\infty, \tilde{\gamma}} \leq \varepsilon(t)\|f\|_{\infty, \gamma}$ for all $f \in \mathcal{B}_\gamma$. Since $\|\cdot\|_{\infty, \gamma} \leq C\|\cdot\|_{\infty, \gamma, \tilde{\gamma}}$ (Proposition 5.2), we obtain

$$|P(t)f - P(t_0)f|_{\tilde{\gamma}} \leq \|P(t)f - P(t_0)f\|_{\infty, \tilde{\gamma}} \leq C\varepsilon(t)\|f\|_{\infty, \gamma, \tilde{\gamma}}.$$

(ii) Let $\rho_K$ be such that $\max\{\vartheta_0^{1/n_0}, r_K\} < \rho_K < 1$, and let $\Gamma$ be the oriented circle $\{|z| = \rho_K\}$ in $\mathbb{C}$. For $t \in K$, we have $r(P(t)) \leq r_K < \rho_K$, thus $P(t)^n = \frac{1}{2i\pi}\int_\Gamma z^n(z - P(t))^{-1}\,dz$. Moreover, the theorem of Keller–Liverani ensures that, for any $t_0 \in K$, there exists an open interval $I$, containing $t_0$, such that $\sup\{\|(z - P(t))^{-1}\|_\gamma, t \in I, |z| = \rho_K\} < +\infty$. By compactness, we get $\sup\{\|(z - P(t))^{-1}\|_\gamma, t \in K, |z| = \rho_K\} < +\infty$. This gives (ii). $\square$

*Proof of assertion 1 of Proposition 9.1.* Let $\gamma$, $\max\{r, s+1\} < \gamma < \gamma_0 - r$. We have $s + 1 < \gamma < \gamma_0$ and $\mathcal{J}^\gamma(0, S) \leq \mathcal{J}^{\gamma_0 - r}(0, S) < +\infty$. Thus Proposition 7.1 applies to $P(t)$ acting on $\mathcal{B}_\gamma$. For convenience, the interval $I_\gamma$ will be denoted by $I$.

LEMMA 9.2. *The maps $v(\cdot)$, $\phi(\cdot)$, and $N(\cdot)$ have derivatives at $t = 0$ as functions with values in $(\mathcal{B}_\gamma, \|\cdot\|_{\infty, \gamma_0})$, $\mathcal{B}'_\gamma$, and $\mathcal{L}(\mathcal{B}_\gamma, \mathcal{B}_{\gamma_0})$, respectively, and there exists a constant $K$ such that we have, for all $n \geq 1$ and all $t \in I$,*

$$\|N(t)^n - N(0)^n\|_{\gamma, \gamma_0} \leq K|t|(\kappa_0')^n.$$

*Moreover, there exists $\gamma_2$, $0 < \gamma_2 < \gamma$, such that $P(\cdot)$ has a derivative at $t = 0$ as an $\mathcal{L}(\mathcal{B}_{\gamma_2}, \mathcal{B}_\gamma)$-valued function.*

PROOF. We have $s + 1 < \gamma \leq \gamma + r < \gamma_0$ and $\mathcal{J}^{\gamma_0 - r}(R, R + S) < +\infty$, so that the condition $\mathcal{V}_1(\gamma, \gamma_0)$ is fulfilled and the assertions upon $v(\cdot)$, $\phi(\cdot)$ and $N(\cdot)$ follow from Proposition 8.1.

Since $r < \gamma < \gamma_0 - r$, there exists $\gamma_2$ such that $0 < \gamma_2 \leq \gamma_2 + r < \gamma \leq \gamma + r < \gamma_0$.

To establish that $P(t)$ has a derivative, we apply Proposition 6.3. Actually, the condition $\mathcal{U}_1(\gamma_2, \gamma)$ holds: we have $\gamma_2 + r < \gamma$, $s + 1 < \gamma$ and $\mathcal{J}^{\gamma_2}(R, R + S) \leq \mathcal{J}^{\gamma_0 - r}(R, R + S) < +\infty$. $\square$

The formula for $\mathbb{E}[e^{itS_n}f(X_n)]$ is obtained by using the basic lemma stated in Section 4.2, the decomposition of $P(t)$ given in Proposition 7.1, and by setting $L(t)f = \langle \phi(t), f \rangle v(t) - \langle \nu, f \rangle 1$.

Under the conditions of (i), we have $\langle \mu, N(0)^n f \rangle = 0$, so that the considered inequality follows from Lemma 9.2.

Inequality (ii) already appears in Corollary 7.2′.

To obtain (iii), it suffices to remark that, since the functions $v(\cdot)$ and $\phi(\cdot)$ have derivatives in $(\mathcal{B}_\gamma, \|\cdot\|_{\infty,\gamma_0})$ and $\mathcal{B}'_\gamma$, there exist constants $C_1$ and $C_2$ such that, for $f \in \mathcal{B}_\gamma$,

$$\|L(t)f\|_{\infty,\gamma_0} \leq |\langle \phi(t), f \rangle| \, \|v(t) - 1\|_{\infty,\gamma_0} \; + \; |\langle \phi(t) - \nu, f \rangle| \, \|1\|_{\infty,\gamma_0}$$
$$\leq C_1|t| \, \|\phi(t)\|_{\infty,\gamma} \|f\|_{\infty,\gamma} \; + \; C_2|t| \, \|f\|_{\infty,\gamma} \|1\|_{\infty,\gamma_0}.$$

It remains to prove the properties of $\lambda(\cdot)$.

From Proposition 7.1, we have $\lambda(0) = 1$ and $\lambda(t)^n = \langle \nu, P(t)^n v(t) \rangle$. Appealing to the invariance of $\nu$, we get $|\lambda(t)|^n \leq \langle \nu, P^n|v(t)| \rangle = \langle \nu, |v(t)| \rangle$. It follows that $|\lambda(t)| \leq 1$.

To prove that $\lambda(\cdot)$ can be expanded to the second order and to identify the terms of its expansion, we proceed as in Lemma IV.4′ of Hennion and Hervé (2001).

LEMMA 9.3. *For $t \in I$, set $p(t) = \langle \phi(t), 1 \rangle$, $\tilde{\nu}(t) = \langle \nu, P(t)1 \rangle$ and $u(t) = P(t)1 - \tilde{\nu}(t)1$. Then $u(0) = 0$, $\langle \nu, u(t) \rangle = 0$, and*

$$\lambda(t) = \frac{1}{p(t)} \langle \phi(t) - \nu, u(t) \rangle + \tilde{\nu}(t).$$

PROOF. The two first equalities are obvious. From the decomposition of Proposition 7.1, we have $P(t)1 = \lambda(t)p(t)v(t) + N(t)1$. As $\langle \phi(t), v(t) \rangle = 1$ and $\phi(t)N(t) = 0$, the formula for $\lambda(t)$ follows from

$$\langle \phi(t), u(t) \rangle = \langle \phi(t), \lambda(t)p(t)v(t) + N(t)1 - \tilde{\nu}(t)1 \rangle = \lambda(t)p(t) - \tilde{\nu}(t)p(t). \quad \square$$

Notice that $\tilde{\nu}(\cdot)$ is the characteristic function of $\xi$ under the distribution $\pi \otimes \nu$, so that the next lemma results from the moment property $\nu(d(\cdot, x_0)^{\gamma_0+1}) < +\infty$.

LEMMA 9.4. *Let $n \in \mathbb{N}^*$. Assume that $\int_G R(g)^n \, d\pi(g) < +\infty$ and that $r \leq \frac{\gamma_0+1}{n}$.*

*Then $\tilde{\nu}(\cdot)$ has continuous derivatives up to order $n$, with $\tilde{\nu}^{(k)}(0) = i^k \int_G \xi(g,x)^k \, d\pi(g) \, d\nu(x)$ for $k = 1, \ldots, n$.*

We can now obtain the second-order Taylor expansion of $\lambda(\cdot)$.

LEMMA 9.5.   $u(\cdot)$ *has a derivative at* $t = 0$ *as a* $\mathcal{B}_\gamma$-*valued function, and we have*

$$\lambda(t) = 1 + imt - \sigma^2 \frac{t^2}{2} + o(t^2),$$

*with*

$$\sigma^2 = (\pi \otimes \nu)(\xi^2) - 2\langle \phi'(0), u'(0) \rangle \geq m^2.$$

PROOF.   By assumption, we have $\int_G R(g)^2 \, d\pi(g) < +\infty$ and $r < \frac{\gamma_0}{2} \leq \frac{\gamma_0 + 1}{2}$, so that $\tilde{\nu}(t) = 1 + imt - (\pi \otimes \nu)(\xi^2)\frac{t^2}{2} + o(t^2)$. From Lemma 9.2, we know that $P(\cdot)1$ and $\langle \nu, P(\cdot)1 \rangle$ have derivatives at $t = 0$ as functions with values in $\mathcal{B}_\gamma$ and $\mathbb{C}$, respectively. Therefore $u(\cdot)$ has a derivative at $t = 0$ as a $\mathcal{B}_\gamma$-valued function.

We get, first in $\mathcal{B}'_\gamma$, $\phi(t) - \nu = \phi(t) - \phi(0) = t\phi'(0) + o(t)$; second in $\mathcal{B}_\gamma$, $u(t) = tu'(0) + o(t)$, and thirdly in $\mathbb{C}$, $p(t) = 1 + O(t)$. Setting $c = 2\langle \phi'(0), u'(0) \rangle$, we have

$$\frac{1}{p(t)} \langle \phi(t) - \nu, u(t) \rangle = \left(1 + O(t)\right)\left(c\frac{t^2}{2} + o(t^2)\right) = c\frac{t^2}{2} + o(t^2).$$

We obtain the Taylor expansion of $\lambda(\cdot)$ by adding the expansion of $\tilde{\nu}$ to the last one.

We now prove that $\sigma^2 \geq m$. Setting $\overline{v(t)}(\cdot) = \overline{v(t)(\cdot)}$, we have $P(-t)\overline{v(t)} = \overline{\lambda(t)} \ \overline{v(t)}$ and, by uniqueness [cf. Proposition 7.1(a)], we get $\lambda(-t) = \overline{\lambda(t)}$. It follows that $\sigma^2 \in \mathbb{R}$. As $1 \geq |\lambda(t)|^2 = 1 - (\sigma^2 - m^2)t^2 + o(t^2)$, we obtain $\sigma^2 - m^2 \geq 0$. Lemma 9.5 is proved.   □

When $m = 0$ and $\sigma^2 > 0$, it follows from the preceding expansion that, for small $|t|$, $|\lambda(\frac{t}{\sigma})| \leq 1 - \frac{t^2}{2} + \frac{t^2}{4} \leq e^{-t^2/4}$, that is, (iv).

*Proof of the assertion* 2 *of Proposition* 9.1.   The claimed inequality follows [cf., e.g., Hennion and Hervé (2001)] from the fact that, under the additional hypotheses in 2, the remainder of the second-order expansion of $\lambda(\cdot)$ can be specified as follows.

PROPOSITION 9.6.   *We have* $\lambda(t) = 1 + imt - \sigma^2 \frac{t^2}{2} + O(t^3)$.

PROOF.   We need the following lemma.

LEMMA 9.7.   *There exists* $0 < \gamma_2 < \gamma_0$ *such that the functions* $\phi(\cdot)$ *and* $P(\cdot)1$ *have a second-order Taylor expansion at* $t = 0$ *as functions with values in* $\mathcal{B}'_{\gamma_2}$ *and in* $\mathcal{B}_{\gamma_2}$, *respectively.*

PROOF. By assumption, we have $\gamma_0 > 3r + \max\{r, s+1\}$; therefore, $4r < \gamma_0$ and $s + 1 + r < \gamma_0 - 2r$. It follows that there exist $\gamma_4$ and $\gamma_2$ such that $0 < \gamma_4 \leq \gamma_4 + 2r < \gamma_2 \leq \gamma_2 + 2r < \gamma_0$ and $s + 1 + r < \gamma_2$.

To establish the assertion on $\phi(\cdot)$, we apply Proposition 8.2. This is possible since the condition $\mathcal{V}_2(\gamma_2, \gamma_0)$ is satisfied; indeed, we have $s + 1 \leq s + 1 + r < \gamma_2 \leq \gamma_2 + 2r < \gamma_0$, and $\mathcal{J}^{\gamma_0 - r}(R, R + S) + \mathcal{J}^{\gamma_2}(R^2, (R + S)R) < +\infty$ because $\gamma_2 < \gamma_0 - 2r$.

Moreover, the condition $\mathcal{U}_2(\gamma_4, \gamma_2)$ is verified: we have $0 < \gamma_4 \leq \gamma_4 + 2r < \gamma_2$, $s + 1 + r < \gamma_2$, and then $\mathcal{J}^{\gamma_4}(R^2, (R + S)R) < +\infty$ since $\gamma_4 < \gamma_0 - 2r$. Proposition 6.3 shows that $P(\cdot)$ has a second-order Taylor expansion at $t = 0$ as an $\mathcal{L}(\mathcal{B}_{\gamma_4}, \mathcal{B}_{\gamma_2})$-valued function, hence the claimed property for $P(\cdot)1$. $\quad\square$

To conclude, we appeal once more to the formula of Lemma 9.3. Since $\int_G R(g)^3 \, d\pi(g) < +\infty$ and $r < \frac{\gamma_0}{4} \leq \frac{\gamma_0 + 1}{3}$, the characteristic function $\tilde{\nu}(\cdot)$ has now three continuous derivatives, so that the remainder of its second-order Taylor expansion is $O(t^3)$. Using the preceding lemma, we have $\phi(t) = \nu + t\phi'(0) + t^2\phi_2 + o(t^2)$ in $\mathcal{B}'_{\gamma_2}$ and $u(t) = tu'(0) + t^2 u_2 + o(t^2)$ in $\mathcal{B}_{\gamma_2}$. Consequently, $\frac{1}{p(t)}\langle \phi(t) - \nu, u(t) \rangle = (1 + O(t))(c\frac{t^2}{2} + O(t^3)) = c\frac{t^2}{2} + O(t^3)$. It follows that the remainder of the second-order expansion of $\lambda(\cdot)$ at $t = 0$ is $O(t^3)$. $\square$

### 9.3. *Proof of Theorem* S'.

PROPOSITION 9.8. *Assume $\mathcal{H}(\gamma_0)$ with $\gamma_0 > r + \max\{r, s+1\}$ and that*

$$\int_G R(g)^2 \, d\pi(g) < +\infty, \qquad \mathcal{J}^{\gamma_0 - r}(R, R + S) < +\infty,$$

*and that $m = 0$.*

(i) *We set $\theta(x) = \int_G \xi(g, x) \, d\pi(g), \; x \in M$. There exists a unique real-valued function $w \in \mathcal{B}_{\gamma_0 - r}$ such that*

$$\langle \nu, w \rangle = 0, \qquad (1 - P)w = \theta,$$

*and we have $\sigma^2 = (\pi \otimes \nu)(\xi(\xi + 2w \circ j))$.*

(ii) *Moreover, suppose that $\gamma_0 > 2r + s + 1$ and that*

$$\mathcal{J}^{\gamma_0 - 2r}(R^2, (R + S)R) < +\infty.$$

*Suppose that the r.v. $Z$ has a distribution $\mu$ which defines an element of $\mathcal{B}'_{\gamma_0}$. Then, for all $n \geq 1$, the characteristic function $\varphi_n(t) = \mathbb{E}[e^{itS_n^Z}] = \langle \mu, P(t)^n 1 \rangle$ has the Taylor expansion $\varphi_n(t) = 1 + a_n t + b_n \frac{t^2}{2} + o_n(t^2)$, with $\sup_{n \geq 1} |b_n + n\sigma^2| < +\infty$.*

Recall that $j$ defines the action of $G$ on $M$.

PROOF OF PROPOSITION 9.8. The hypothesis $\int_G R(g)^2 \, d\pi(g) < +\infty$ implies that $\theta$ is well defined.

*Proof of assertion* (i).    To begin, we state the differential properties that we shall use.

Let $\gamma$ be such that $\max\{r, s+1\} < \gamma < \gamma_0 - r$. Then there exist $\gamma_2, \gamma'$ such that $0 < \gamma_2 \leq \gamma_2 + r < \gamma' < \gamma \leq \gamma + r < \gamma_0$ and $s + 1 < \gamma'$. It is easily checked that we have the following properties and their consequences:

1. $\mathcal{V}_1(\gamma, \gamma_0)$; therefore $\phi(\cdot)$ has a derivative at $t = 0$ as a $\mathcal{B}'_\gamma$-valued function (Proposition 8.1).
2. $\mathcal{U}_1(\gamma_2, \gamma')$; therefore $P(\cdot)$ has the derivative $L_1$ at $t = 0$ as an $\mathcal{L}(\mathcal{B}_{\gamma_2}, \mathcal{B}_{\gamma'})$-valued function (Proposition 6.3).
3. $\mathcal{U}_1(\gamma', \gamma_0)$; therefore $P(\cdot)$ has the derivative $L_1$ at $t = 0$ as an $\mathcal{L}(\mathcal{B}_{\gamma'}, \mathcal{B}_{\gamma_0})$-valued function (Proposition 6.3).
4. $\mathcal{U}_0(\gamma', \gamma)$; therefore $P(\cdot)$ is continuous at $t = 0$ as an $\mathcal{L}(\mathcal{B}_{\gamma'}, \mathcal{B}_\gamma)$-valued function (Proposition 6.2).

Lemma 9.4 asserts that $\tilde{\nu}(\cdot)$ has a continuous derivative, with $\tilde{\nu}'(0) = im = 0$. The property 2 above ensures that $u(\cdot)$ has a derivative at $t = 0$ as a $\mathcal{B}_{\gamma'}$-valued function, and that $u'(0)(x) = L_1 \mathbf{1}(x) - \tilde{\nu}'(0) = i \int_G \xi(g, x) \, d\pi(g) = i\theta(x)$; thus $u'(0) = i\theta$. It follows that $\theta \in \mathcal{B}_{\gamma'}$. Since $\langle \nu, \theta \rangle = im = 0$, Theorem 5.5 shows that there exists a unique $w \in \mathcal{B}_{\gamma'}$ such that $\langle \nu, w \rangle = 0$ and $(1 - P)w = \theta$, and that $w$ is the sum in $\mathcal{B}_{\gamma'}$ of the series $\sum_{n \geq 0} P^n \theta$. As $\theta$ is real valued, so is $w$. At last, since $\gamma' < \gamma_0 - r$, we have $w \in \mathcal{B}_{\gamma_0 - r}$.

On the basis of the formula of Lemma 9.5, we get

$$\sigma^2 = (\pi \otimes \nu)(\xi^2) - 2i\langle \phi'(0), \theta \rangle = (\pi \otimes \nu)(\xi^2) - 2i\langle \phi'(0), (1 - P)w \rangle.$$

The following lemma allows us to conclude.

LEMMA 9.9.  *We have*

$$\langle \phi'(0), (1 - P)w \rangle = i(\pi \otimes \nu)(\xi \, w \circ j).$$

PROOF.    It is known that, for small $|t|$, $(\lambda(t) - P(t))^* \phi(t) = \phi(t)(\lambda(t) - P(t)) = 0$. Hence, setting $S(t) = \lambda(t) - P(t)$, we have

$$\left( \frac{\phi(t) - \phi(0)}{t} \right) S(t)w + \phi(0) \left( \frac{S(t)w - S(0)w}{t} \right) = \frac{\phi(t)S(t)w - \phi(0)S(0)w}{t} = 0.$$

Observe that $\lambda(t)$ has a derivative at $t = 0$ because the conditions of point (1) in Proposition 9.1 hold. Therefore, since $w \in \mathcal{B}_{\gamma'}$ and $\phi(0) = \nu \in \mathcal{B}'_{\gamma_0}$, the above properties (1) and (4), and then (3) enable us to pass to the limit in the equality. We get

$$\phi'(0)(S(0)w) + \nu(S'(0)w) = 0,$$

or else $\phi'(0)(1 - P)w = \nu[(L_1 - \lambda'(0))w] = \nu(L_1 w) = i \int_M \int_G \xi(g, x) \times w(gx) \, d\pi(g) \, d\nu(x)$. $\square$

*Proof of assertion* (ii). We know that $\lambda(t) = 1 - \sigma^2 \frac{t^2}{2} + o(t^2)$. Otherwise, since $s + 1 < \gamma_0 - 2r$, there exists $\eta$ such that $s + 1 < \eta \leq \eta + 2r < \gamma_0$. Since $\mathcal{J}^{\gamma_0 - r}(R, R + S) < +\infty$ and $\mathcal{J}^\eta(R^2, (R + S)R) < +\infty$ ( because $\eta < \gamma_0 - 2r$), the condition $\mathcal{V}_2(\eta, \gamma_0)$ holds. Consequently, Proposition 8.2 applies. It follows that $v(\cdot)$, $\phi(\cdot)$, $N(\cdot)$ have second-order Taylor expansions at $t = 0$ as functions with values in $(\mathcal{B}_\eta, N_{\infty, \gamma_0})$, in $\mathcal{B}'_\eta$, and in $\mathcal{L}(\mathcal{B}_\eta, \mathcal{B}_{\gamma_0})$, respectively. We get, for all $n \geq 1$,

$$\langle \phi(t), 1 \rangle \langle \mu, v(t) \rangle = 1 + tB + \frac{t^2}{2} C + o(t^2) \qquad (A, B \in \mathbb{C}),$$

$$\langle \mu, N(t)^n 1 \rangle = \langle \mu, N(0)^n 1 \rangle + t \langle \mu, N_{1,n} 1 \rangle + \frac{t^2}{2} \langle \mu, N_{2,n} 1 \rangle + o_n(t^2).$$

Since $N(0)1 = 0$ and $\varphi_n(t) = \langle \mu, P(t)^n 1 \rangle = \lambda(t)^n \langle \phi(t), 1 \rangle \langle \mu, v(t) \rangle + \langle \mu, N(t)^n 1 \rangle$, with $\lambda(t)^n = 1 - n\sigma^2 \frac{t^2}{2} + o_n(t^2)$, the coefficient $b_n$ of $\frac{t^2}{2}$ in the Taylor expansion of $\varphi_n$ is $C - n\sigma^2 + \langle \mu, N_{2,n} 1 \rangle$. This enables us to conclude because $\sup_{n \geq 1} \|N_{2,n}\|_{\eta, \gamma_0} < +\infty$.   $\square$

*End of the proof of Theorem* S$'$.

*Proof of* (ii). Let us prove that $\mathbb{E}[(S_n^Z)^2] < +\infty$. Actually, since $2r \leq \gamma_0$, we have, for $k \geq 1$,

$$\mathbb{E}[\xi(Y_k, Z_{k-1})^2] \leq \mathbb{E}[R(Y_k)^2] \mathbb{E}[\psi(Z_{k-1})] = \int_G R^2 \, d\pi \int_M P^{k-1} \psi \, d\mu,$$

with $\psi(x) = (1 + d(x, x_0))^{\gamma_0}$. Since $\psi \in \mathcal{B}_{\gamma_0 - 1} \subset \mathcal{B}_{\gamma_0}$, $P \in \mathcal{L}(\mathcal{B}_{\gamma_0})$, and $\mu \in \mathcal{B}'_{\gamma_0}$, we get $\mathbb{E}[\xi(Y_k, Z_{k-1})^2] < +\infty$; hence the claimed property. The function $\varphi_n(\cdot)$ has therefore a second-order derivative at $t = 0$ and $\varphi_n''(0) = -\mathbb{E}[(S_n^Z)^2]$. With the help of Proposition 9.8(ii), we obtain $\mathbb{E}[(S_n^Z)^2] = -b_n$, hence $\sigma^2 = \lim_n \frac{1}{n} \mathbb{E}[(S_n^Z)^2]$.

*Proof of* (i). The method of the proof of Theorem IV.7 of Hennion and Hervé (2001) applies here to the transition probability $Q$ introduced in Section 4.2, yet we give below an adaptation of this method only using $P$. Set

$$\tilde{\xi} = \xi + w \circ j,$$

where $w$ is the function in Proposition 9.8 (it can be checked that $\tilde{\xi} - Q\tilde{\xi} = \xi$). Recall that $\sigma^2 = (\pi \otimes \nu)[\xi(\xi + 2w \circ j)]$ (Proposition 9.8). From the equality $\xi^2 + 2\xi \, w \circ j = (\xi + w \circ j)^2 - (w \circ j)^2 = \tilde{\xi}^2 - (w \circ j)^2$, we get

$$\sigma^2 = (\pi \otimes \nu)(\tilde{\xi}^2 - (w \circ j)^2).$$

Assume that $\nu(w^2) < +\infty$. Then, using the invariance of $\nu$, we can write

$$\sigma^2 = (\pi \otimes \nu)(\tilde{\xi}^2) - \nu(w^2) = \int_M d\nu(x) \int_G (\tilde{\xi}(g, x)^2 - w(x)^2) \, d\pi(g).$$

But

$$\int_G \tilde{\xi}(g,x)\,d\pi(g) = \int_G \xi(g,x)\,d\pi(g) + Pw(x) = \theta(x) + (w(x) - \theta(x)) = w(x),$$

so that

$$\sigma^2 = \int_M d\nu(x) \int_G (\tilde{\xi}(g,x) - w(x))^2\,d\pi(g)$$

$$= \int_M d\nu(x) \int_G (\xi(g,x) + w(gx) - w(x))^2\,d\pi(g).$$

If $\sigma^2 = 0$, we therefore get $\xi(g,x) = w(x) - w(gx)$  $\pi \otimes \nu$ a.e.

To complete the proof of Theorem S′, it now suffices to show that the hypothesis $\sigma^2 = 0$ implies $\nu(w^2) < +\infty$. We know that, for all $x \in M$, $w(x) = \sum_{n \geq 0} P^n \theta(x)$. Since $P^n \theta(x) = \mathbb{E}[\theta(R_n x)] = \mathbb{E}[\int \xi(g, R_n x)\,d\pi(g)] = \mathbb{E}[\xi(Y_{n+1}, R_n x)]$, we have, for all $x \in M$, $w(x) = \lim_n \mathbb{E}[S_n^x]$.

Assume that $Z$ has the distribution $\nu$ and that $\sigma^2 = 0$. Then the point (ii) of Proposition 9.8 and the fact that $b_n = -\mathbb{E}_\mu[(S_n^Z)^2]$ show that $\sup_n \mathbb{E}[(S_n^Z)^2] = \vartheta < +\infty$. From the inequalities $\int \mathbb{E}[S_n^x]^2\,d\nu(x) \leq \int \mathbb{E}[(S_n^x)^2]\,d\nu(x) = \mathbb{E}[(S_n^Z)^2]$ and Fatou's Lemma, we deduce that $\nu(w^2) \leq \vartheta$.

EXAMPLE (Study of $\sigma^2$ for sequences of type $(u(Y_n)\chi(Z_{n-1}))_n$).   Suppose that the function $\xi$ is of the form $\xi(g,x) = u(g)\chi(x)$, where $u$ is a nonzero real valued measurable function on $G$ and $\chi$ is a real-valued locally Lipschitz function on $M$ satisfying $|\chi(x) - \chi(y)| \leq Cd(x,y)(1 + d(x,x_0) + d(y,x_0))^s$. Observe that Condition RS holds with $r = s+1$ and $R(s) = S(g) = C|u(g)|$.

In this context, the next statement, based on both Theorem S′ and 5.5, gives a simple sufficient condition for $\sigma^2 > 0$.

PROPOSITION S″.   *Suppose that the conditions of Theorem* S *hold* [*with $r = s+1$ and $R(s) = S(g) = |u(g)|$*], *that $\int_G u(g)\,d\pi(g) = 0$, and that $\chi(x) \neq 0$ for some $x$ in the support $\Sigma_\nu$ of the $P$-invariant measure $\nu$. Then $\sigma^2 > 0$.*

PROOF.   Observe that $m = \pi \otimes \nu(\xi) = 0$. By Theorem S′, we shall get $\sigma^2 > 0$ if we prove that there is no real-valued function $\tilde{\chi}_1$ in $\mathcal{B}_{\gamma_0 - r}$ such that, for all $x \in \Sigma_\nu$, we have $\xi(g,x) = u(g)\chi(x) = \tilde{\xi}_1(x) - \tilde{\xi}_1(gx)$ $\pi$-a.e.

Let $\tilde{\xi}_1$ be such a function. Then, by integrating the above equality with respect to the measure $\pi$, we get $\tilde{\xi}_1(x) = \int_G \tilde{\xi}_1(gx)\,d\pi(g) = (P\tilde{\xi}_1)(x)$ for all $x \in \Sigma_\nu$. Since $\Sigma_\nu$ is an absorbing set [for all $x \in \Sigma_\nu$, we have $P(x, \Sigma_\nu) = 1$], this can be rewritten as $\tilde{\xi}_{1|\Sigma_\nu} = P_{\Sigma_\nu}(\tilde{\xi}_{1|\Sigma_\nu})$, where $P_{\Sigma_\nu}$ denotes the kernel induced by $P$ on $\Sigma_\nu$. From $\mathrm{Ker}(P-1) = \mathbb{C} \cdot 1$ (Theorem 5.5), it can be easily proved that the functions of $\mathcal{B}_{\gamma_0 - r}$ whose restriction on $\Sigma_\nu$ are $P_{\Sigma_\nu}$-invariant are constant on $\Sigma_\nu$. It follows that $\tilde{\xi}_{1|\Sigma_\nu}$ is constant; thus, for all $x \in \Sigma_\nu$, $u(g)\chi(x) = 0$ $\pi$-a.e. This is impossible.   □

Institut Mathematique de Rennes
Université de Rennes 1
UMR-CNRS 6625
Campus de Beaulieu
35042 Rennes Cedex
France
e-mail: hubert.hennion@univ-rennes1.fr

Institut National des Sciences
   Appliquées de Rennes
20 Avenue des Buttes de
   Coësmes CS 14 315
35043 Rennes Cedex
France
e-mail: loic.herve@insa-rennes.fr